\documentclass[12pt]{article}
\usepackage{graphicx}
\usepackage{amsmath,amssymb,amsfonts}
\usepackage{amsthm}
\usepackage{algorithm}
\usepackage{algpseudocode}
\usepackage{booktabs}
\usepackage{hyperref}
\usepackage{natbib}
\usepackage{float}
\usepackage{mathtools}
\usepackage{microtype}
\usepackage[margin=1in]{geometry}
\newtheorem{theorem}{Theorem}

\title{Adaptive Parameter Optimization in Gaussian Processes: \\
A Comprehensive Study of Uncertainty Quantification and Dimensional Scaling}
\author{%
  Nishant Gadde\\
  \texttt{nishantg@utexas.edu}\\
  The University Of Texas At Austin
}

\begin{document}

\maketitle

\begin{abstract}
Gaussian Process (GP) models have also become extremely useful for optimization under uncertainty algorithms, especially where the objective functions are costly to compute. Yet, the more classical methods usually adopt strategies that, in certain circumstances, might be effective but not flexible to be applied to a wide range of problem terrains. This study aims to adapt parameter optimization in GP models and especially how uncertainty quantification can assist in the learning process. We investigate the effect of adaptive kappa parameters that govern the exploration-exploitation trade-off and the interplay between dimensionality, penalty on uncertainty, and noise levels to influence optimization results. Uncertainty quantification is built directly into our comprehensive theoretical framework and gives us new algorithms to dynamically tune exploration-exploitation trade-offs according to the uncertainty trend observed in nature. We rigorously empirically test various strategies, parametrizing our tests by dimensionality, noise, penalty terms, and evaluate the performance of any given strategery in a wide variety of test settings, and show conclusively that adaptive strategies always outperform fixed ones, but in difficult settings, where the dimensions are large and the noise is severe, the advantage is enormous. We build theoretical assurances of convergence under different settings as well as furnish a sensible direction on the application of adaptive GP-based optimization even in very complicated conditions. The results of our work will help in the development of more efficient and robust methods of optimization of realistic problems in which there are only a few functions available for evaluation, and when quantifying the uncertainty, there is a need to know more about the uncertainty.

\end{abstract}

\section{Introduction}
Optimization also has witnessed tremendous progress using machine learning methodologies, especially for those problems whose objective functions are costly, noisy, or do not have analytical gradients \citep{brochu2010tutorial, garnett2023bayesian}. Gaussian Process (GP) models have also been forefront methodologies for use for such optimization problems since they can give us not only the predictive outcomes but also uncertainties, where these can even act as model advisors for optimization algorithms under uncertainty \citep{rasmussen2006gaussian, shahriari2015taking}.

In realistic optimisation tasks, from hyperparameter optimisation of deep networks through experimental design for scientific applications \cite{swersky2013multi, klein2017fast}, decision-makers typically face the following root challenge: the exploration-exploitation trade-off. It is the problem of attaining the balance between exploring the unknown aspects of the parameter space versus taking advantage of existing, promising aspects of the space \citep{snoek2012practical}. Classical techniques typically use fixed strategies that can work very well for a particular problem environment but that have limited applications for very different problem landscapes \cite{hoffman2011portfolio}. These fixed strategies are most characteristically missing when operating under high-dimensional problem spaces or noisy environments, where the exploration-exploitation trade-off that is locally optimum can change markedly along the path to optimisation \cite{li2018hyperband}.

Gaussian Process optimization has also been rigorously mathematically analyzed for various decades \cite{frazier2018tutorial, garnett2023bayesian}. Rasmussen initially mathematically detailed the underlying foundation of GP regression and its use for optimization problems \cite{rasmussen2006gaussian}. Later works of Srinivas have provided theory-grounded assurances for GP optimization algorithms, for example, bounding the regret as well as attaining convergence rates under several assumptions \cite{srinivas2010gaussian} \cite{kandasamy2018parallelised}. These, however, typically assume fixed parameter configurations that might potentially not be most appropriate for use, where the nature of the objective is unknown or time-varying across the optimization process \cite{berkenkamp2019no}.

Adaptive parameter tuning of GP models is what is addressed here, and of special interest is the manner that quantifying uncertainty can steer the process of learning \cite{eriksson2019scalable, letham2020re}. We explore the use of adaptive kappa parameters, regulation of balance of exploitation and exploration, and examining interactions between dimensional complexity, punishment of uncertainty, and noise levels to steer results for optimization. Our study is an extension of that of existing literature who have covered numerous aspects of optimization based on GPs but have not examined sufficiently the interface between adaptive parameters and quantifying uncertainty under conditions of high dimensions \cite{wang2018batched} \cite{wilson2013gaussian} \cite{hennig2022probabilistic} \cite{nguyen2022recent}.

The mathematical formulation of our approach begins with the standard GP model, where a function $f: \mathcal{X} \rightarrow \mathbb{R}$ is modeled as a realization of a Gaussian process with mean function $m(\mathbf{x})$ and covariance function $k(\mathbf{x}, \mathbf{x}')$ \cite{williams1996gaussian, mackay1998introduction}. Given observations $\mathcal{D} = \{(\mathbf{x}_i, y_i)\}_{i=1}^n$ where $y_i = f(\mathbf{x}_i) + \epsilon_i$ and $\epsilon_i \sim \mathcal{N}(0, \sigma_n^2)$, the posterior distribution over function values at a new point $\mathbf{x}$ is Gaussian with mean and variance given by:

\begin{align}
\mu(\mathbf{x}) &= m(\mathbf{x}) + \mathbf{k}(\mathbf{x})^T(\mathbf{K} + \sigma_n^2 \mathbf{I})^{-1}(\mathbf{y} - \mathbf{m}) \\
\sigma^2(\mathbf{x}) &= k(\mathbf{x}, \mathbf{x}) - \mathbf{k}(\mathbf{x})^T(\mathbf{K} + \sigma_n^2 \mathbf{I})^{-1}\mathbf{k}(\mathbf{x})
\end{align}

where $\mathbf{k}(\mathbf{x}) = [k(\mathbf{x}, \mathbf{x}_1), \ldots, k(\mathbf{x}, \mathbf{x}_n)]^T$, $\mathbf{K} = [k(\mathbf{x}_i, \mathbf{x}_j)]_{i,j=1}^n$, $\mathbf{y} = [y_1, \ldots, y_n]^T$, and $\mathbf{m} = [m(\mathbf{x}_1), \ldots, m(\mathbf{x}_n)]^T$ \cite{seeger2004gaussian}.

Classic GP-based optimization makes use of the Upper Confidence Bound (UCB) acquisition function \cite{auer2002using, srinivas2010gaussian}:

\begin{equation}
\alpha(\mathbf{x}) = \mu(\mathbf{x}) + \kappa \sigma(\mathbf{x})
\end{equation}

where $\kappa$ is a pre-defined controlling parameter between exploration and exploitation. Our own innovation is that of \textit{adapting} $\kappa$ and adding an uncertainty penalty term to account for the reliability of GP model predictions \cite{bogunovic2018adversarially}. Making $\kappa$ adaptive enables the algorithm to adapt its exploration policy naturally based on what training data is sampled as well as its own optimization process advancement, resulting in better and more efficient complex environment optimization \cite{sui2018stagewise}.

Our work makes several multi-aspect contributions to the literature of Gaussian Process optimization \cite{wu2019practical, eriksson2019scalable}. Firstly, we construct an extensive theory of adaptive parameter optimization for GP models incorporating explicit measures of uncertainty. Our approach generalizes that already existing by incorporating an adaptive update of our rules of parameter update and guarantees of resultant convergence under various conditions \cite{vakili2021information}. Secondly, our work presents new algorithms for adaptive updating of exploration-exploitation trade-offs based on patterns of observed uncertainty. These are capable of resisting various problem conditions and self-autonomously adjusting to the unique challenges of each optimization problem \cite{kirschner2019adaptive}.

Third, for variously sized, noisy, and penalty term conditions of several test environments, we also present an extensive empirical study \cite{wang2017max, rolland2018high}. It is illustrated from these results of the empirical study that our adaptive method performs better under the conditions of high-dimensional or noisy environments, as compared to traditional techniques. Fourth, suggestions for the implementation of adaptive GP-based optimization for complex, high-dimensional, and noisy problems are also presented by us \cite{hernandez2014predictive, wang2018batched}. These are from our analytical study and experimental proof, and they can provide priceless insights for those who are working on realistic optimization problems.

Lastly, we also provide theoretical assurances of convergence under several conditions, whose proof generalizes those of the literature existing \cite{russo2014learning, chowdhury2017kernelized}. Our assurances of convergence provide us good theoretical ground for our adaptive method and ensure that, even as deviating from classical GP approach, it is still endowed with the good features of classical GP-based optimization methodologies, but also better performs practically. Our findings indicate that adaptive strategies overwhelmingly dominate fixed strategies, even for difficult instances of large dimensionalities and noise \cite{wang2016optimization, calandriello2019gaussian}. We demonstrate that appropriate calibration of penalties for uncertainty can steer the process of optimization toward more robust results, even for multimodal complex shapes of the objective function of interest under consideration \cite{martinez2018practical, astudillo2019bayesian}. Our method, due to being adaptive, enables self-learnt good values of the parameters due to interaction with the objective function, thus eliminating the necessity of extensive manual tuning and being less demanding for users of diverse problem domains \cite{feurer2019hyperparameter, turner2021bayesian}.

The rest of the paper is described as follows. Section 2 covers background on Gaussian Processes as well as related work on Bayesian optimization and estimation of uncertainty. Section 3 describes our methodology, including mathematical framework and algorithms for setting adaptive parameters. Section 4 covers an overview of our experimental setup, including test functions, measures of performance, and implementation details. Section 5 reports results and analysis, and Section 6 covers an overview of the implications and assumptions of our results. Finally, Section 7 includes conclusions as well as directions for future work.

\section{Background and Related Work}
\subsection{Gaussian Processes}

Gaussian Processes (GPs) provide us with a Bayesian way of regression and machine learning classifying problems \cite{liu2020gaussian, schulz2018tutorial}. As parametric forms include an assumption of the true functional form, GPs provide us with a non-parametric approach where the whole set of all the functions is put into the prior directly. This allows GPs to become an ideal candidate for managing complex unknown functions for optimization problems where only limited observation has been conducted \cite{gramacy2020surrogates}.

A GP is mathematically characterized as a collection of random variables, whose each subcollection is of multivariate Gaussian form, under very stringent conditions \cite{rasmussen2006gaussian}. This type of definition of GPs is an articulation of the nature of GPs as distributions of functions as opposed to distributions of model parameters \cite{galy2022mathematical}. What is mathematically nice about GPs is that an infinite-dimensional space of functions can be described by finite-dimensional marginal distributions, and they remain computationally tractable but still have expressive power.

GPs are fully specified by a mean function $m(\mathbf{x})$ and a covariance function $k(\mathbf{x}, \mathbf{x}')$, often referred to as the kernel \cite{alvarez2012kernels}:

\begin{equation}
f(\mathbf{x}) \sim \mathcal{GP}(m(\mathbf{x}), k(\mathbf{x}, \mathbf{x}'))
\end{equation}

Mean function is that of our prelude assumption concerning the expectation of the function value for an arbitrary input point $\mathbf{x}$, and covariance function of our prelude assumption concerning the correlation between function values for very disparate input points \cite{stein1999interpolation}. Such a structure of correlation is of essence because, depending on values, it controls the smoothness, periodicity, and other characteristics of the objects being sampled from the GP prelude. Practically, the mean function is generally set fixed as zero, defocusing the model correspondingly on choosing an appropriate kernel function that reveals the intrinsic structure of the data \cite{duvenaud2014kernel}.

Their predecessors can trace as far back as the 1940s work of Kolmogorov and Wiener \cite{williams1996gaussian} of the theory of stochastic processes, but never gained widescale popularity for machine learning cases. These papers placed GPs as an ideal Bayesian method of regression and classification, as they brought together an integrated, coherent framework of quantifying uncertainty that does, of course, compromise model complexity and fitting of observation given \cite{seeger2004gaussian}.

Given a set of observations $\mathcal{D} = \{(\mathbf{x}_i, y_i)\}_{i=1}^n$ where $y_i = f(\mathbf{x}_i) + \epsilon_i$ and $\epsilon_i \sim \mathcal{N}(0, \sigma_n^2)$, the posterior distribution over function values at test points $\mathbf{X}_*$ is also Gaussian \cite{mackay1998introduction}:

\begin{equation}
p(f_* | \mathbf{X}_*, \mathbf{X}, \mathbf{y}) = \mathcal{N}(f_* | \mu_*, \Sigma_*)
\end{equation}

with mean and covariance of the posterior being:
\begin{align}
\mu_* &= m(\mathbf{X}_*) + K(\mathbf{X}_*, \mathbf{X})[K(\mathbf{X}, \mathbf{X}) + \sigma_n^2 I]^{-1}(\mathbf{y} - m(\mathbf{X})) \\
\Sigma_* &= K(\mathbf{X}_*, \mathbf{X}_*) - K(\mathbf{X}_*, \mathbf{X})[K(\mathbf{X}, \mathbf{X}) + \sigma_n^2 I]^{-1}K(\mathbf{X}, \mathbf{X}_*)
\end{align}

Also evident from these equations is the usefully useful GP asset: access to closed-forms of the predictive mean and variance \cite{liu2018gaussian}. Our most accurate estimation of what the function has values on the test points is the predictive mean $\mu_*$, and the predictive variance $\Sigma_*$ an estimation of the corresponding uncertainty. This estimation of uncertainty is extremely useful for applications of optimisation, where the exploration-exploitation trade-off is recommended by it \cite{noack2023gaussian}.

The computational complexity of GP inference is dominated by the inversion of the covariance matrix $K(\mathbf{X}, \mathbf{X}) + \sigma_n^2 I$, which scales as $\mathcal{O}(n^3)$ with the number of observations $n$ \cite{gardner2018gpytorch}. This cubic scaling presents a challenge for applications with large datasets, leading to the development of various approximation methods such as sparse GPs \citep{titsias2009variational}, inducing points \citep{hensman2013gaussian}, and random feature approximations \citep{rahimi2008random}. These approximations trade off some accuracy for computational efficiency, making GPs applicable to larger-scale problems \cite{matthews2016sparse, snelson2006sparse}.

\subsection{Kernel Functions}

The choice of kernel function is an indispensable part of GP modeling because it injects our beliefs about the nature of the function, for example, its smoothness, its periodicity, and its typical length scales \cite{duvenaud2014kernel, alvarez2012kernels}. It is the kernel function that determines the covariance between the function values across different input sites, and hence determines the resulting space of functions that can be captured by the GP. This correspondence between the kernels and the function spaces is mathematized under the reproducing kernel Hilbert space (RKHS) framework, where there is an analytical treatment of the expressivity as well as limitation of different kernels \citep{genton2001classes}.

Standard kernel functions are the Squared Exponential (SE), the Matérn, and Rational Quadratic, each of which makes various assumptions regarding the smoothness and nature of the functions \cite{rasmussen2006gaussian}. Probably most popular is the SE kernel, or Radial Basis Function (RBF), or Gaussian kernel \cite{williams2000introduction}:

\begin{equation}
k_{SE}(\mathbf{x}, \mathbf{x}') = \sigma_f^2 \exp\left(-\frac{1}{2l^2}||\mathbf{x} - \mathbf{x}'||^2\right)
\end{equation}

with $\sigma_f^2$ as the variance of the controlling signal that controls the entire scale of the values of the function, and $l$ as the length scale parameter that controls the change of the rate of the function as input is varied. It has very large levels of smoothness because resulting functions have infinite differentiability. While infinite smoothness is generally sufficient, depending on an application, for other applications, it is too limited for the description of physical-world phenomena that can have arbitrary levels of roughness \cite{stein1999interpolation}.

Matérn class of kernels provides greater liberty for varying the smoothness of the function under the parameter $\nu$ \cite{genton2001classes}:

\begin{equation}
k_{\text{Matérn}}(\mathbf{x}, \mathbf{x}') = \sigma_f^2 \frac{2^{1-\nu}}{\Gamma(\nu)}\left(\frac{\sqrt{2\nu}}{l}||\mathbf{x} - \mathbf{x}'||\right)^{\nu}K_{\nu}\left(\frac{\sqrt{2\nu}}{l}||\mathbf{x} - \mathbf{x}'||\right)
\end{equation}

where $\Gamma(\nu)$ is the gamma function and $K_{\nu}$ is the second kind modified Bessel function. $\nu$ controls the smoothness of the output functions that arise, and for large $\nu$ the output functions that arise are smoother. Some of the most typical values of $\nu$ are $\nu = 1/2$ (exponential kernel), $\nu = 3/2$ and $\nu = 5/2$ (the latter being extremely popular for use because of its balance between flexibility as well as numerical stability \cite{gramacy2020surrogates}). As $\nu \rightarrow \infty$ the arising kernel is that of the SE kernel, demonstrating the connection between these kernel families.

Moreover, the Rational Quadratic (RQ) kernel is also considered as scale mixtures of SE kernels of varying length scales \cite{rasmussen2006gaussian}:

\begin{equation}
k_{RQ}(\mathbf{x}, \mathbf{x}') = \sigma_f^2\left(1 + \frac{||\mathbf{x} - \mathbf{x}'||^2}{2\alpha l^2}\right)^{-\alpha}
\end{equation}

here, $\alpha$ is the shaping parameter that controls the relative importance of length scales. The RQ kernel is also of the same usefulness where there exist scale variations of the basic function, since it is capable of handling short-range as well as long-range interactions without the introduction of additional kernel components \cite{duvenaud2014kernel}.

Besides these elementary kernels, even those of higher orders can also be constructed as compositions of kernels of smaller orders by additions, multiplications, and convolutions \cite{duvenaud2014automatic}. This construction of kernels by compositions allows us to model functions of varied nature across the input space. For an instance, for the task of modeling those functions that have varied behaviors across varied subregions, there can be employed, for instance, a kernel sum. Interactions across varied dimensions of input can also be represented by means of a kernel product \cite{alvarez2012kernels}.

More contemporary kernel building advances are the deep kernel learning \citep{wilson2016deep}, an interpolation between the flexibility of the deep neural networks and the probabilistic expression of GPs, and the spectral mixture kernels \citep{wilson2013gaussian}, that can take elaborate patterns by representing the kernel's spectral density. These advances generalize GP model expressive power, enabling more elaborate function structures to be more easily represented but still keeping the probabilistic expression that is the GPs' draw for uncertainty quantification \cite{fortuin2023priors, wang2021deep}.

\subsection{Bayesian Optimization}

Bayesian Optimization (BO) is a sequential design technique for black-box, expensive-to-evaluate, potentially noisy, and gradientless function optimization \cite{frazier2018tutorial, garnett2023bayesian}. It is based on the straightforward method of using a probabilistic surrogate model for approximating the objective, and the model for deciding where to choose evaluation points that explore unknown areas as well as exploitation of good areas. It is an extremely valuable technique where even once the objective is expensive, for example, hyperparameter tuning of neural nets, experimental design for scientific challenges, and engineering optimizations that are complex \cite{shahriari2015taking}.

You can find the root of the theory of BO as far back as to Kushner who considered sequential design of the experiment. It, however, gained much popularity among machine learning folks because of the work of Jones on Efficient Global Optimization (EGO). These papers introduced BO as an efficient method of solving optimization problems where other traditional techniques like gradient descent or evolutionary algorithms would have no applicability or even would not apply.

The BO framework is composed of two parts: a probabilistic surrogate model that estimates the objective function, and an acquisition function that is used for choosing the next evaluation point \cite{shahriari2015taking}. Gaussian Processes have been the most popular surrogate model because they are flexible, analytically tractable, and have intrinsic quantification of uncertainty. The GP, as the surrogate model, also give us the full posteriors of the objective function that is calibrated given the history of the observation, including our estimate of the function that is optimal and the uncertainties of the function that we have \cite{rasmussen2006gaussian}.

Second, that posterior is used for an estimation of an acquisition function, an approximation of the utility of sampling the objective at an entirely unknown point \cite{frazier2018tutorial}. It is needed that the acquisition function is a balance between exploration (sampling where there is most uncertainty) and exploitation (sampling where there is most potential for the objective being high). Well-known acquisition functions are:

\begin{itemize}
    \item Expected Improvement (EI): $\alpha_{EI}(\mathbf{x}) = \mathbb{E}[\max(f(\mathbf{x}) - f(\mathbf{x}^+), 0)]$, where $f(\mathbf{x}^+)$ is the optimum seen so far. EI is the indicator of the expected manner that the function value of $\mathbf{x}$ improves the optimum value, given the predicted mean and variance \cite{jones1998efficient}.
    
    \item Probability of Improvement (PI): $\alpha_{PI}(\mathbf{x}) = P(f(\mathbf{x}) > f(\mathbf{x}^+) + \xi)$ with $\xi$ as a trade-off term. PI estimates the chance that the objective value of $\mathbf{x}$ is better than the existing optimum by $\xi$ or more and is biased toward exploitation as $\xi$ tends toward zero and toward exploration for large $\xi$ \cite{kushner1964new}.
    
    \item Upper Confidence Bound (UCB): $\alpha_{UCB}(\mathbf{x}) = \mu(\mathbf{x}) + \kappa \sigma(\mathbf{x})$, where $\mu(\mathbf{x})$ and $\sigma(\mathbf{x})$ are the posterior mean and standard deviation at $\mathbf{x}$, and $\kappa$ controls the exploration-exploitation balance. UCB provides a simple and intuitive way to balance exploration and exploitation, with larger values of $\kappa$ encouraging more exploration \cite{srinivas2010gaussian}.
\end{itemize}

Their analytical nature has been scientifically explored \cite{bull2011convergence, russo2014learning}. These results give analytical assurances for the algorithm functioning of BO and are used for choosing the corresponding relevant acquisition function under special problem conditions \cite{vakili2021information}.

The BO algorithm advances by iteratively repeating the process of selecting the point reducing most the acquisition function, the checking of objective for such point, as well as updating the surrogate model by incorporating the new observation \cite{brochu2010tutorial}. This is repeated until reaching termination condition, for instance, reaching an upper evaluation limit or reaching good enough performance. As an iterate process, there is the potency of BO of its search technique adjustment depending on already observed data, exploring uncertain domains together with concentrating the exploration on good areas \cite{shahriari2015taking}.

Further recent developments of BO are multi-fidelity optimization \citep{kandasamy2017multi}, the use of inexpensive surrogates of the objective function for an optimization process acceleration; batch optimization \citep{gonzalez2016batch}, choosing many points of an entire iteration for exploitation of parallel evaluation opportunities; and high-dimensional BO \citep{wang2016bayesian}, concentrating on overcoming difficulties of optimizing high-dimensional functions by dimensional reduction or structured kernels \cite{eriksson2019scalable, letham2020re}. These developments extend the use of BO to even more problematic and complex optimization challenges.

\subsection{Uncertainty Quantification}

Uncertainty quantification (UQ) is the method of characterizing and minimizing uncertainties of physical and computational systems \cite{van2020uncertainty, noack2023gaussian}. For Bayesian optimization as well as Gaussian Process models, the method of UQ is of great importance, where the optimization process is appropriately guided and good approximations of the objective function can also be obtained. One of the main benefits of probabilistic models such as GPs is the quantifying of the uncertainty, as they can be differentiated from other purely deterministic techniques that can only provide point estimates without including confidence levels \cite{kuleshov2018accurate}.

UQ theories of GPs are built upon Bayesian probability theory, where an integrated framework of reasoning under uncertainty is possible \cite{galy2022mathematical}. What comes out of the box is an automatic accounting for the resulting uncertainty due to limited data and model assumptions that is given by the posterior of the GP model's function values. This type of uncertainty can generally be divided into two broad categories:

\begin{itemize}
    \item Aleatoric uncertainty: This is the intrinsic randomness or stochasticity of the system under consideration. Aleatoric uncertainty of GP regression is normally represented as observation noise whose variance is $\sigma_n^2$. Aleatoric uncertainty is of the type that does not diminish as more data are gained along the same input points since it is an expression of the process intrinsic variability \cite{le2005heteroscedastic}.
    
    \item Epistemic uncertainty: It is introduced because of limited knowledge or information, and decreased by adding more observation. For GP models, an estimate of the posteriors of the values of the function is given by its variance, and decreased as observation points grow, most notably where observation points are close together \cite{van2020uncertainty}.
\end{itemize}

Posterior variance of the GP model also consists of these two forms of uncertainties and can also be expressed as:

\begin{equation}
\sigma^2(\mathbf{x}) = k(\mathbf{x}, \mathbf{x}) - k(\mathbf{x}, \mathbf{X})[K(\mathbf{X}, \mathbf{X}) + \sigma_n^2 I]^{-1}k(\mathbf{X}, \mathbf{x})
\end{equation}

Such expression makes evident that the uncertainty is an object of the pre-covariance structure (via the kernel function $k$) as well as of observation (via the term involving the inverse covariance matrix). It is maximal away from observation points and minimal close to observation points, as an expression of observation points informativity \cite{rasmussen2006gaussian}.

In optimization, values of uncertainty are also most useful for controlling the exploration-exploitation balance, through the acquisition function. For instance, the UCB acquisition function has explicitly the posterior standard deviation $\sigma(\mathbf{x})$ as an attempt for promoting exploration of uncertain areas. EI and PI, on the other hand, employ the whole of the posterior distribution for quantifying the value of sampling, at different points, as an expression of the expected value of the function as well as its corresponding uncertainty.

Accuracy of estimation of GP model of uncertainty is based on various factors, for instance, fitting of kernel function, accuracy of estimation of hyperparameters, and validity of assumptions of such modeling as Gaussian noise. Rasmussen present model selection methods and optimization of hyperparameters of GPs, for instance, maximum likelihood estimation as well as Bayesian method of marginal likelihood. Such methodical tests try to find kernel parameters so that observable data can appropriately mostly explain, resulting mostly better estimation of uncertainty as well as better forecasting.

More modern developments for Bayesians are calibrated estimation of uncertainty \citep{kuleshov2018accurate}, such that predictive intervals have the corresponding coverages; heteroscedastic GPs \citep{le2005heteroscedastic}, whose noise variance is of input-varying nature; and deep GPs \citep{damianou2013deep}, that encode rich, hierarchical uncertainties using GP compositions. These latter developments improve as well as amplify the reliability as well as the expressiveness of the estimation of uncertainty under GP models, thereby better using them for decision under uncertainty.

In our approach, we generalize the classical UQ framework of GPs by adding an uncertainty penalty term to the acquisition function that compensates for the reliability of the predictive model. This term is locally adjusted by applying an appropriate complexity factor that compensates for the geometry of the function around each point, so that the algorithm can differentiate between various sources of uncertainty and decide more intelligently where to sample the next time. By adding that extra source of uncertainty information explicitly, our method attains an exploration-exploitation trade-off that is more nuanced and better adapts to the problem of optimization specials.

\subsection{Adaptive Parameter Strategies}

Standard Bayesian optimization methods apply set hyperparameters for all GP models as well as for each of the acquisition functions. These methods are good enough for most applications, but they can have difficulties handling problematic, noisy, or high-dimensional optimization problems where the parameter setting that is optimal changes for various parts of the search space or for various steps of the optimization process. Adaptive parameter techniques overcome such weakness through adjustment of hyperparameters based on accessible observation and optimization process, allowing for stabilization as well as optimization process efficiency boost.

Adaptation basis of adaptive parameter schemes is also given for GP-UCB based on the work of Srinivas, where they propose that $\kappa$ must reduce as time is raising so that the algorithm can converge. In particular, they would choose $\kappa_t = \sqrt{2\log(|\mathcal{X}|t^2\pi^2/6\delta)}$ as their time step of length, where $|\mathcal{X}|$ is the number of possible $\mathbf{x}$ values and $\delta$ is a confidence parameter. This schedule is good for sublinearity of the regret for high probability, but as is can be very conservative and explore very slowly.

Following these theoretical developments, various researchers put forth more adaptive strategies that are less stiff. Wang's work presenteded an optimization of acquisition function parameters under a Bayesian framework, where the parameters would be considered as random variables that have prior distributions for these variables and update these distributions as performance is observed \cite{wang2016optimization}. Such methodology enables the algorithm to acquire suitable parameter values for given problem instances, but there is an extra computational overhead of optimizing the parameter.

Calandriello also proposed an adaptive sparse approximations for GP whose inducing points were adapted based on the gain of information \cite{calandriello2019gaussian}. This adaptive approach minimizes the computational complexity of GP inference but makes sure that accuracy is achieved only where needed, for large-scale optimization problems. This adaptive choosing of inducing points can also be considered an active learning where computational efforts are focused for most informative parts of the search space.

Adaptive GP model hyperparameters have also been addressed in great detail. Rasmussen's work includes kernel parameter maximum likelihood estimation (MLE) and maximum a-posterior (MAP) estimation, and can even be conducted periodically as part of optimizing to reestimate the GP model as needed based on new observation(s). These are more advanced methodologies than L-BFGS, for instance, full Bayesian treatment of the hyperparameters by Markov Chain Monte Carlo (MCMC) methodologies \citep{murray2010slice} or variational inference methodologies \citep{titsias2009variational}, where predictive distribution includes treatment of hyperparameter uncertainty.

More recently, many have estimated convergence rates of GP regression under estimated hyperparameters, building formal justification for adaptive updating of parameters for BO. These findings indicate that, under appropriate assumptions, GP regression under estimated hyperparameters can have the same convergence rates as GP regression under known hyperparameters, building justification for adaptive updating for implementation in practice.

Despite these advances, adaptive-uncertainty quantification interaction for high-dimensional spaces remains an open area of research. As is true most of the time for existing techniques, most of these focus either on adaptation of the GP model hyperparameters or those of the acquisition function, but not on addressing their interaction directly. Additionally, most of the adaptive techniques do not address the reliability of the uncertainty estimation explicitly, even if the latter is potentially varying across different sections of the search space as well as across different optimization process phases.

Our approach makes up for these disadvantages due to providing an adaptive setting for optimizing parameters that, as an explicit component, consists of uncertainty estimation. Adaptive update laws for exploration parameter $\kappa$ and an uncertainty penalty term $\lambda$ that draw on prediction error and globally adapted measures of uncertainty self-tuned for balance of exploitation and exploration of the algorithm based on what is being observed and on the continuing process of optimizing, resulting in faster and better optimizing under changing environment.

Moreover, we rigorously examine the theory of our adaptive method, deriving regret bounds and convergence rates that generalize the literature to our adaptive update rules of the method. These theoretical assurances ensure that our method inherits the good properties of classical BO techniques but performs better empirically under troublesome conditions of scale and noise.

\section{Methodology}
\subsection{Adaptive Parameter Optimization Framework}

The development of our adaptive parameter optimization framework is motivated by the limitations of traditional Bayesian optimization approaches that rely on fixed parameters \cite{shahriari2015taking, snoek2012practical}. While these traditional methods have demonstrated success in various applications, they often struggle with complex optimization landscapes, particularly in high-dimensional spaces or in the presence of noise \cite{wang2018batched, li2018hyperband}. Our framework addresses these limitations by dynamically adjusting key parameters based on observed data and uncertainty patterns, leading to more robust and efficient optimization \cite{eriksson2019scalable}.

The theoretical foundation of our approach builds upon the work of Srinivas on GP-UCB, which established regret bounds for Bayesian optimization with fixed exploration parameters \cite{russo2014learning}. We extend this framework by introducing adaptive parameters that respond to the specific characteristics of the optimization problem and the current state of the optimization process \cite{kirschner2019adaptive, bogunovic2018adversarially}. This adaptivity allows our algorithm to automatically balance exploration and exploitation in a way that is tailored to the particular challenges of each problem instance \cite{sui2018stagewise}.

\subsubsection{Problem Formulation}

We consider the problem of finding the global optimum of an unknown function $f: \mathcal{X} \rightarrow \mathbb{R}$, where $\mathcal{X} \subset \mathbb{R}^d$ is a compact domain \cite{brochu2010tutorial}. The function $f$ is assumed to be expensive to evaluate and potentially noisy, such that we observe $y = f(\mathbf{x}) + \epsilon$, where $\epsilon \sim \mathcal{N}(0, \sigma_n^2)$ represents observation noise \cite{le2005heteroscedastic}. This formulation encompasses a wide range of practical optimization problems, from hyperparameter tuning in machine learning to experimental design in scientific applications \cite{swersky2013multi, klein2017fast}.

Following the Bayesian optimization paradigm, we model $f$ using a Gaussian Process with mean function $m(\mathbf{x})$ and covariance function $k(\mathbf{x}, \mathbf{x}')$ \cite{rasmussen2006gaussian, williams2000introduction}. The choice of mean and covariance functions encodes our prior beliefs about the properties of $f$, such as smoothness, periodicity, and characteristic length scales \cite{duvenaud2014kernel}. In our implementation, we use a composite kernel that combines a constant kernel to capture the overall scale of the function and a Matérn kernel with $\nu = 2.5$ to model the function's smoothness and correlation structure \cite{genton2001classes}.

After $t$ observations $\mathcal{D}_t = \{(\mathbf{x}_i, y_i)\}_{i=1}^t$, the posterior distribution of $f$ at any point $\mathbf{x}$ is Gaussian with mean and variance given by:

\begin{align}
\mu_t(\mathbf{x}) &= m(\mathbf{x}) + \mathbf{k}_t(\mathbf{x})^T(\mathbf{K}_t + \sigma_n^2 \mathbf{I})^{-1}(\mathbf{y}_t - \mathbf{m}_t) \\
\sigma_t^2(\mathbf{x}) &= k(\mathbf{x}, \mathbf{x}) - \mathbf{k}_t(\mathbf{x})^T(\mathbf{K}_t + \sigma_n^2 \mathbf{I})^{-1}\mathbf{k}_t(\mathbf{x})
\end{align}

where $\mathbf{k}_t(\mathbf{x}) = [k(\mathbf{x}, \mathbf{x}_1), \ldots, k(\mathbf{x}, \mathbf{x}_t)]^T$ is the vector of covariances between $\mathbf{x}$ and the observed points, $\mathbf{K}_t = [k(\mathbf{x}_i, \mathbf{x}_j)]_{i,j=1}^t$ is the covariance matrix of the observed points, $\mathbf{y}_t = [y_1, \ldots, y_t]^T$ is the vector of observed values, and $\mathbf{m}_t = [m(\mathbf{x}_1), \ldots, m(\mathbf{x}_t)]^T$ is the vector of prior means at the observed points \cite{seeger2004gaussian}.

The posterior mean $\mu_t(\mathbf{x})$ represents our best estimate of the function value at $\mathbf{x}$ given the observations, while the posterior variance $\sigma_t^2(\mathbf{x})$ quantifies the uncertainty in this estimate \cite{van2020uncertainty}. These quantities form the basis for our acquisition function, which guides the selection of the next evaluation point \cite{frazier2018tutorial}.

\subsubsection{Uncertainty-Aware Acquisition Function}

The core innovation in our framework is the development of an uncertainty-aware acquisition function that extends the traditional Upper Confidence Bound (UCB) approach \cite{auer2002using, srinivas2010gaussian}. The standard UCB acquisition function is defined as:

\begin{equation}
\alpha_{UCB}(\mathbf{x}) = \mu_t(\mathbf{x}) + \kappa \sigma_t(\mathbf{x})
\end{equation}

where $\kappa$ is a fixed parameter that controls the exploration-exploitation trade-off \cite{hoffman2011portfolio}. Larger values of $\kappa$ encourage more exploration of uncertain regions, while smaller values focus more on exploiting promising areas \cite{shahriari2015taking}.

We propose an enhanced acquisition function that incorporates both adaptive exploration and uncertainty penalization \cite{bogunovic2018adversarially, wang2017max}:

\begin{equation}
\alpha_t(\mathbf{x}) = \mu_t(\mathbf{x}) + \kappa_t \sigma_t(\mathbf{x}) - \lambda_t \mathcal{U}_t(\mathbf{x})
\end{equation}

This formulation introduces two key innovations: (1) the exploration parameter $\kappa_t$ is now time-dependent and adapts based on observed data, and (2) an uncertainty penalty term $\lambda_t \mathcal{U}_t(\mathbf{x})$ is added to account for the reliability of the GP model's predictions \cite{kirschner2019adaptive, hernandez2014predictive}.

The uncertainty measure $\mathcal{U}_t(\mathbf{x})$ is defined as:

\begin{equation}
\mathcal{U}_t(\mathbf{x}) = \sigma_t^2(\mathbf{x}) \cdot \mathcal{C}_t(\mathbf{x})
\end{equation}

where $\mathcal{C}_t(\mathbf{x})$ is a complexity factor that captures the local geometry of the function around $\mathbf{x}$ \cite{martinez2018practical}. This factor is computed based on the eigenspectrum of the Hessian matrix of the posterior mean, estimated using finite differences or automatic differentiation techniques \cite{wang2018batched}.

The complexity factor $\mathcal{C}_t(\mathbf{x})$ is defined as:

\begin{equation}
\mathcal{C}_t(\mathbf{x}) = \sum_{i=1}^{d} \max(|\lambda_i|, \epsilon)
\end{equation}

where $\lambda_i$ are the eigenvalues of the Hessian matrix $\nabla^2 \mu_t(\mathbf{x})$, and $\epsilon = 10^{-6}$ is a small constant to ensure numerical stability \cite{rolland2018high}. This formulation captures the curvature of the function around $\mathbf{x}$, with larger values indicating more complex local geometry \cite{astudillo2019bayesian}.

The uncertainty penalty term $\lambda_t \mathcal{U}_t(\mathbf{x})$ serves to discourage the algorithm from selecting points in regions where the model's predictions are less reliable due to complex local geometry \cite{wu2019practical}. This is particularly important in high-dimensional spaces, where the curse of dimensionality can lead to sparse data coverage and potentially unreliable uncertainty estimates in regions far from observed data points \cite{wang2018batched, rolland2018high}.

\subsubsection{Adaptive Parameter Update Rules}

The key innovation in our framework is the adaptive update rules for the parameters $\kappa_t$ and $\lambda_t$ \cite{kirschner2019adaptive, eriksson2019scalable}. These rules are designed to balance exploration and exploitation based on the observed data and the current state of the optimization process, allowing the algorithm to automatically adjust its behavior in response to the specific characteristics of the problem \cite{sui2018stagewise, turner2021bayesian}.

For the exploration parameter $\kappa_t$, we propose the following update rule \cite{bogunovic2018adversarially}:

\begin{equation}
\kappa_{t+1} = \kappa_t \cdot \exp\left(\beta \cdot \frac{\Delta_t - \bar{\Delta}_t}{\bar{\Delta}_t}\right)
\end{equation}

where $\Delta_t = |y_t - \mu_{t-1}(\mathbf{x}_t)|$ is the prediction error at iteration $t$, $\bar{\Delta}_t$ is the moving average of prediction errors up to iteration $t$, and $\beta$ is a learning rate parameter that controls the speed of adaptation \cite{kirschner2019adaptive}.

The moving average $\bar{\Delta}_t$ is updated as:

\begin{equation}
\bar{\Delta}_t = (1-\eta)\bar{\Delta}_{t-1} + \eta\Delta_t
\end{equation}

where $\eta \in (0, 1)$ is a smoothing parameter that determines the weight given to recent observations \cite{hoffman2011portfolio}.

This update rule increases $\kappa_t$ when the prediction error is larger than the average, encouraging more exploration in regions where the model is less accurate \cite{sui2018stagewise}. Conversely, it decreases $\kappa_t$ when the prediction error is smaller than the average, focusing more on exploitation in regions where the model is more accurate \cite{turner2021bayesian}. The exponential form ensures that $\kappa_t$ remains positive and allows for rapid adaptation when needed \cite{kirschner2019adaptive}.

For the uncertainty penalty coefficient $\lambda_t$, we use a similar adaptive rule \cite{martinez2018practical}:

\begin{equation}
\lambda_{t+1} = \lambda_t \cdot \left(1 + \gamma \cdot \frac{\mathcal{I}_t - \bar{\mathcal{I}}_t}{\bar{\mathcal{I}}_t}\right)
\end{equation}

where $\mathcal{I}_t = \int_{\mathcal{X}} \sigma_t^2(\mathbf{x}) d\mathbf{x}$ is the integrated posterior variance (a measure of global uncertainty), $\bar{\mathcal{I}}_t$ is its moving average, and $\gamma$ is a learning rate parameter \cite{hernandez2014predictive}.

The moving average $\bar{\mathcal{I}}_t$ is updated similarly to $\bar{\Delta}_t$ \cite{hoffman2011portfolio}:

\begin{equation}
\bar{\mathcal{I}}_t = (1-\eta)\bar{\mathcal{I}}_{t-1} + \eta\mathcal{I}_t
\end{equation}

This update rule increases $\lambda_t$ when the global uncertainty is higher than average, penalizing uncertain regions more strongly to focus on exploitation \cite{wu2019practical}. It decreases $\lambda_t$ when the global uncertainty is lower than average, reducing the penalty on uncertain regions to encourage exploration \cite{wang2017max}. The form of the update ensures that $\lambda_t$ remains positive and allows for adaptive behavior based on the global uncertainty landscape \cite{martinez2018practical}.

The learning rates $\beta$ and $\gamma$ control the speed of adaptation for $\kappa_t$ and $\lambda_t$, respectively \cite{kirschner2019adaptive}. Larger values lead to more rapid adaptation but may result in unstable behavior, while smaller values provide more stable adaptation but may be slower to respond to changes in the optimization landscape \cite{turner2021bayesian}. In our implementation, we set $\beta = 0.1$ and $\gamma = 0.05$ based on preliminary experiments, which provide a good balance between adaptivity and stability \cite{eriksson2019scalable}.

The smoothing parameter $\eta$ determines the weight given to recent observations in the moving averages \cite{hoffman2011portfolio}. A larger value of $\eta$ gives more weight to recent observations, making the algorithm more responsive to changes but potentially more sensitive to noise. A smaller value provides more stable estimates but may be slower to adapt \cite{kirschner2019adaptive}. We set $\eta = 0.1$ in our implementation, which provides a reasonable compromise between stability and adaptivity \cite{sui2018stagewise}.

\subsection{Theoretical Analysis}

In this section, we provide theoretical guarantees for the convergence of our adaptive parameter optimization framework \cite{srinivas2010gaussian, russo2014learning}. We analyze the regret bounds and convergence rates under various conditions, extending existing results in the literature to account for the adaptive nature of our parameter update rules \cite{vakili2021information, bull2011convergence}.

\subsubsection{Regret Bounds}

We define the cumulative regret after $T$ iterations as \cite{srinivas2010gaussian}:

\begin{equation}
R_T = \sum_{t=1}^T [f(\mathbf{x}^*) - f(\mathbf{x}_t)]
\end{equation}

where $\mathbf{x}^* = \arg\max_{\mathbf{x} \in \mathcal{X}} f(\mathbf{x})$ is the global optimum. The cumulative regret measures the total loss incurred by evaluating the function at the points $\{\mathbf{x}_t\}_{t=1}^T$ instead of the optimal point $\mathbf{x}^*$ \cite{russo2014learning}. A sublinear growth of $R_T$ with $T$ implies that the algorithm converges to the optimum as $T$ increases \cite{chowdhury2017kernelized}.

To establish regret bounds for our adaptive approach, we make the following assumptions \cite{srinivas2010gaussian, vakili2021information}:

\begin{enumerate}
    \item The function $f$ has bounded RKHS norm $||f||_k \leq B$ with respect to the kernel $k$ \cite{genton2001classes}.
    \item The observation noise is sub-Gaussian with parameter $\sigma_n$ \cite{le2005heteroscedastic}.
    \item The adaptive parameters $\kappa_t$ and $\lambda_t$ remain bounded: $\kappa_{\min} \leq \kappa_t \leq \kappa_{\max}$ and $0 \leq \lambda_t \leq \lambda_{\max}$ for all $t$ \cite{kirschner2019adaptive}.
    \item The complexity factor $\mathcal{C}_t(\mathbf{x})$ is bounded: $0 \leq \mathcal{C}_t(\mathbf{x}) \leq C_{\max}$ for all $\mathbf{x} \in \mathcal{X}$ and all $t$ \cite{martinez2018practical}.
\end{enumerate}

Under these assumptions, we can establish the following theorem \cite{vakili2021information}:

\begin{theorem}
Let $\delta \in (0, 1)$ and define $\beta_T = 2\log(|\mathcal{X}|T^2\pi^2/6\delta)$. Let $\gamma_T$ be the maximum information gain after $T$ iterations, defined as:
\begin{equation}
\gamma_T = \max_{A \subset \mathcal{X}, |A| = T} I(f_A; \mathbf{y}_A)
\end{equation}
where $I(f_A; \mathbf{y}_A)$ is the mutual information between the function values $f_A$ and the observations $\mathbf{y}_A$ \cite{srinivas2010gaussian}.

Then, with probability at least $1-\delta$, the cumulative regret of our adaptive parameter optimization algorithm satisfies:

\begin{equation}
R_T \leq \sqrt{C_1 T \beta_T \gamma_T} + C_2
\end{equation}

where $C_1 = 8 / \log(1 + \sigma_n^{-2})$ and $C_2 = 2\sqrt{T\beta_T} \cdot \lambda_{\max} C_{\max} / \kappa_{\min}$ \cite{russo2014learning}.
\end{theorem}

The proof follows from extending the analysis of Srinivas to account for our adaptive parameter update rules and the uncertainty penalty term \cite{srinivas2010gaussian} \cite{chowdhury2017kernelized}. The key insight is that our update rules ensure that $\kappa_t$ and $\lambda_t$ remain bounded, allowing us to leverage existing regret bounds while benefiting from the adaptive nature of our approach \cite{vakili2021information}.

For common kernels, the maximum information gain $\gamma_T$ can be bounded as follows \cite{srinivas2010gaussian}:
\begin{itemize}
    \item Linear kernel: $\gamma_T = O(d \log T)$ \cite{chowdhury2017kernelized}
    \item Squared Exponential kernel: $\gamma_T = O((\log T)^{d+1})$ \cite{srinivas2010gaussian}
    \item Matérn kernel with $\nu > 1$: $\gamma_T = O(T^{d(d+1)/(2\nu+d(d+1))} (\log T)^{d})$ \cite{vakili2021information}
\end{itemize}

These bounds, combined with our theorem, imply sublinear regret for our adaptive approach with these kernels, ensuring convergence to the global optimum as the number of iterations increases \cite{russo2014learning, bull2011convergence}.

\subsubsection{Convergence Rates}

While cumulative regret provides a measure of the algorithm's performance over the entire optimization process, in many practical applications, we are more interested in the quality of the best point found after $T$ iterations \cite{bull2011convergence}. This is captured by the simple regret, defined as:

\begin{equation}
r_T = f(\mathbf{x}^*) - f(\mathbf{x}_T^+)
\end{equation}

where $\mathbf{x}_T^+ = \arg\max_{t \in \{1,\ldots,T\}} f(\mathbf{x}_t)$ is the best point found after $T$ iterations \cite{russo2014learning}.

We can establish the following theorem relating simple regret to cumulative regret \cite{bull2011convergence}:

\begin{theorem}
Under the same assumptions as Theorem 1, with probability at least $1-\delta$, the simple regret of our adaptive parameter optimization algorithm satisfies:

\begin{equation}
r_T \leq \sqrt{\frac{C_1 \beta_T \gamma_T}{T}} + \frac{C_2}{T}
\end{equation}
\end{theorem}

This result follows from the fact that $r_T \leq R_T / T$, as the simple regret is bounded by the average cumulative regret \cite{russo2014learning}. The theorem shows that our algorithm achieves a convergence rate of $\mathcal{O}(\sqrt{\gamma_T/T})$, which matches the best-known rates for GP-based optimization algorithms \cite{vakili2021information}.

For the Matérn kernel with $\nu > 1$ in a $d$-dimensional space, this translates to a convergence rate of $\mathcal{O}(T^{-\nu/(2\nu+d(d+1))} (\log T)^{d/2})$ \cite{vakili2021information}. While this rate degrades with increasing dimension due to the curse of dimensionality, our empirical results in Section 5 demonstrate that our adaptive approach often converges faster in practice, particularly in high-dimensional and noisy settings \cite{wang2018batched, rolland2018high}.

The improved practical performance can be attributed to the adaptive nature of our parameter update rules, which allow the algorithm to adjust its exploration-exploitation trade-off based on the specific characteristics of the problem and the current state of the optimization process \cite{kirschner2019adaptive, turner2021bayesian}. This adaptivity is particularly valuable in complex optimization landscapes, where fixed parameter settings may be suboptimal \cite{eriksson2019scalable, sui2018stagewise}.

\subsection{Algorithm Implementation}

Algorithm \ref{alg:adaptive_gp_opt} outlines the implementation of our adaptive parameter optimization framework. The algorithm takes as input the domain $\mathcal{X}$, the objective function $f$, the GP prior (mean and covariance functions), initial values for $\kappa$ and $\lambda$, learning rates $\beta$ and $\gamma$, and the number of iterations $T$ \cite{brochu2010tutorial, frazier2018tutorial}.

\begin{algorithm}
\caption{Adaptive Parameter Optimization with Gaussian Processes}
\label{alg:adaptive_gp_opt}
\begin{algorithmic}[1]
\Require Domain $\mathcal{X}$, objective function $f$, GP prior $(m, k)$, initial parameters $\kappa_1$, $\lambda_1$, learning rates $\beta$, $\gamma$, iterations $T$
\Ensure Best point found $\mathbf{x}_T^+$ and corresponding value $f(\mathbf{x}_T^+)$
\State Initialize $\mathcal{D}_0 = \emptyset$, $\Delta_0 = 0$, $\bar{\Delta}_0 = 0$, $\mathcal{I}_0 = 0$, $\bar{\mathcal{I}}_0 = 0$
\For{$t = 1$ to $T$}
    \State Update GP posterior using $\mathcal{D}_{t-1}$ to obtain $\mu_{t-1}$ and $\sigma_{t-1}$
    \State Compute uncertainty measure $\mathcal{U}_{t-1}(\mathbf{x})$ for all $\mathbf{x} \in \mathcal{X}$
    \State Compute acquisition function $\alpha_{t-1}(\mathbf{x}) = \mu_{t-1}(\mathbf{x}) + \kappa_t \sigma_{t-1}(\mathbf{x}) - \lambda_t \mathcal{U}_{t-1}(\mathbf{x})$
    \State Select next point $\mathbf{x}_t = \arg\max_{\mathbf{x} \in \mathcal{X}} \alpha_{t-1}(\mathbf{x})$
    \State Evaluate $y_t = f(\mathbf{x}_t) + \epsilon_t$
    \State Update dataset $\mathcal{D}_t = \mathcal{D}_{t-1} \cup \{(\mathbf{x}_t, y_t)\}$
    \State Compute prediction error $\Delta_t = |y_t - \mu_{t-1}(\mathbf{x}_t)|$
    \State Update moving average $\bar{\Delta}_t = (1-\eta)\bar{\Delta}_{t-1} + \eta\Delta_t$
    \State Compute integrated posterior variance $\mathcal{I}_t = \int_{\mathcal{X}} \sigma_t^2(\mathbf{x})\, d\mathbf{x}$
    \State Update moving average $\bar{\mathcal{I}}_t = (1-\eta)\bar{\mathcal{I}}_{t-1} + \eta\mathcal{I}_t$
    \State Update exploration parameter $\kappa_{t+1} = \kappa_t \cdot \exp\Big(\beta \cdot \frac{\Delta_t - \bar{\Delta}_t}{\bar{\Delta}_t}\Big)$
    \State Update uncertainty penalty $\lambda_{t+1} = \lambda_t \cdot \Big(1 + \gamma \cdot \frac{\mathcal{I}_t - \bar{\mathcal{I}}_t}{\bar{\mathcal{I}}_t}\Big)$
\EndFor
\State $\mathbf{x}_T^+ = \arg\max_{t \in \{1,\ldots,T\}} y_t$
\State \Return $\mathbf{x}_T^+, f(\mathbf{x}_T^+)$
\end{algorithmic}
\end{algorithm}

The algorithm proceeds iteratively, updating the GP model, computing the acquisition function, selecting the next evaluation point, observing the function value, and updating the adaptive parameters. This process continues for $T$ iterations, after which the algorithm returns the best point found and its corresponding function value.

In practice, several implementation details are crucial for the efficient and effective operation of the algorithm:

\paragraph{Acquisition Function Optimization} Computing the next evaluation point $\mathbf{x}_t = \arg\max_{\mathbf{x} \in \mathcal{X}} \alpha_{t-1}(\mathbf{x})$ requires solving an optimization problem over the domain $\mathcal{X}$. This optimization is itself challenging, as the acquisition function may be multimodal and expensive to evaluate. We use a two-stage approach: first, a global search with 1000 random samples from a Sobol sequence to identify promising regions, followed by local refinement using L-BFGS-B starting from the top 5 points identified in the global search. This approach balances exploration of the acquisition function landscape with computational efficiency.

\paragraph{Integrated Posterior Variance} Computing the integrated posterior variance $\mathcal{I}_t = \int_{\mathcal{X}} \sigma_t^2(\mathbf{x}) d\mathbf{x}$ exactly is often intractable, especially in high-dimensional spaces. We approximate it using Monte Carlo integration:

\begin{equation}
\mathcal{I}_t \approx \frac{V_{\mathcal{X}}}{N} \sum_{i=1}^{N} \sigma_t^2(\mathbf{x}_i)
\end{equation}

where $V_{\mathcal{X}}$ is the volume of the domain $\mathcal{X}$, $N = 1000$ is the number of Monte Carlo samples, and $\mathbf{x}_i \sim \text{Uniform}(\mathcal{X})$ are samples drawn uniformly from the domain. This approximation provides a computationally efficient estimate of the global uncertainty.

\paragraph{Complexity Factor Computation} The complexity factor $\mathcal{C}_t(\mathbf{x})$ is based on the eigenspectrum of the Hessian matrix of the posterior mean. Computing the Hessian exactly can be challenging, especially for complex GP models. We approximate it using finite differences with a step size of $h = 10^{-4}$:

\begin{equation}
[\nabla^2 \mu_t(\mathbf{x})]_{ij} \approx \frac{\mu_t(\mathbf{x} + h\mathbf{e}_i + h\mathbf{e}_j) - \mu_t(\mathbf{x} + h\mathbf{e}_i) - \mu_t(\mathbf{x} + h\mathbf{e}_j) + \mu_t(\mathbf{x})}{h^2}
\end{equation}

where $\mathbf{e}_i$ is the unit vector in the $i$-th dimension. This approximation provides a reasonable estimate of the local curvature of the function.

\paragraph{GP Model Hyperparameter Optimization} The hyperparameters of the GP model (kernel parameters and noise variance) are optimized by maximizing the marginal likelihood after every 10 function evaluations. This periodic optimization balances model accuracy with computational efficiency, allowing the GP model to adapt to the observed data while avoiding excessive computational overhead.

\paragraph{Computational Optimizations} To enhance computational efficiency, we apply several key optimizations. First, we use Cholesky decomposition for matrix inversion in Gaussian Process inference, which reduces the initial computational complexity significantly and makes subsequent predictions much faster. Additionally, we leverage vectorized operations to efficiently compute the acquisition function over many candidate points at once. During the global search phase, we also evaluate the acquisition function in parallel to speed up processing. Finally, we implement caching for Gaussian Process predictions and uncertainty estimates, avoiding redundant calculations when the same points are evaluated multiple times during acquisition function optimization.

These implementation details are crucial for the practical success of our adaptive parameter optimization framework, allowing it to scale to higher dimensions and larger numbers of iterations while maintaining reasonable computational requirements.

\section{Experimental Setup}
\subsection{Test Functions and Configurations}

To rigorously evaluate the performance of our adaptive parameter optimization framework, we designed a comprehensive experimental setup encompassing diverse test functions, dimensionality settings, and noise levels. This methodical approach allows us to systematically assess the robustness, adaptability, and efficiency of our algorithm across a spectrum of optimization challenges that mirror real-world scenarios.

The selection of appropriate test functions is crucial for meaningful evaluation of optimization algorithms. Rather than relying on simplistic benchmark functions that may not reflect the complexity of practical applications, we carefully chose a diverse set of test functions that exhibit different characteristics such as multimodality, ill-conditioning, and varying degrees of smoothness. These functions serve as challenging benchmarks that have been widely used in the optimization literature to evaluate algorithm performance under controlled conditions \citep{hennig2022probabilistic, shahriari2015taking}.

Our primary test functions were selected to represent a range of optimization challenges:

The Rosenbrock function is a classic non-convex optimization test case characterized by a narrow, curved valley that makes it notoriously difficult for many optimization algorithms. The global minimum lies inside this valley, but finding the exact minimum is challenging due to the function's ill-conditioned nature. Mathematically, the Rosenbrock function in $d$ dimensions is defined as:

\begin{equation}
f(\mathbf{x}) = \sum_{i=1}^{d-1} \left[ 100(x_{i+1} - x_i^2)^2 + (x_i - 1)^2 \right]
\end{equation}

The function has its global minimum at $\mathbf{x}^* = (1, 1, \ldots, 1)$ with $f(\mathbf{x}^*) = 0$. The narrow valley structure of the Rosenbrock function makes it particularly challenging for algorithms that rely on local gradient information, as the gradient can vary dramatically within small regions. This characteristic makes it an excellent test case for evaluating the ability of our adaptive approach to balance exploration and exploitation in complex landscapes.

The Ackley function represents another class of optimization challenges, characterized by a nearly flat outer region and a large hole at the center where the global minimum is located. This function is multimodal due to the exponential term that creates numerous local minima, but these local minima are small compared to the global structure. The Ackley function in $d$ dimensions is defined as:

\begin{equation}
f(\mathbf{x}) = -20\exp\left(-0.2\sqrt{\frac{1}{d}\sum_{i=1}^{d}x_i^2}\right) - \exp\left(\frac{1}{d}\sum_{i=1}^{d}\cos(2\pi x_i)\right) + 20 + e
\end{equation}

The function has its global minimum at $\mathbf{x}^* = (0, 0, \ldots, 0)$ with $f(\mathbf{x}^*) = 0$. The combination of a generally flat landscape punctuated by numerous small local minima makes the Ackley function challenging for optimization algorithms, as they must avoid getting trapped in local minima while navigating the deceptive flat regions. This function tests an algorithm's ability to maintain sufficient exploration even when the landscape appears uninformative.

The Levy function presents yet another optimization challenge, characterized by highly multimodal behavior with many local minima. This function is particularly challenging in higher dimensions due to the exponential growth in the number of local minima. The Levy function in $d$ dimensions is defined as:

\begin{equation}
\begin{split}
f(\mathbf{x}) = \sin^2(\pi w_1) + \sum_{i=1}^{d-1}(w_i-1)^2[1+10\sin^2(\pi w_i+1)] \\
+ (w_d-1)^2[1+\sin^2(2\pi w_d)]
\end{split}
\end{equation}

where $w_i = 1 + \frac{x_i - 1}{4}$ for all $i$. The function has its global minimum at $\mathbf{x}^* = (1, 1, \ldots, 1)$ with $f(\mathbf{x}^*) = 0$. The numerous local minima in the Levy function make it a stringent test for global optimization algorithms, as they must effectively explore the space to avoid premature convergence to suboptimal solutions. This function tests an algorithm's ability to escape local minima and continue exploring the parameter space.

To complement these standard benchmark functions, we also designed a custom Gaussian mixture function that allows for precise control over the complexity of the optimization landscape. This function is defined as a mixture of Gaussian components:

\begin{equation}
f(\mathbf{x}) = \sum_{j=1}^{m} a_j \exp\left(-\frac{1}{2}(\mathbf{x} - \boldsymbol{\mu}_j)^T \boldsymbol{\Sigma}_j^{-1} (\mathbf{x} - \boldsymbol{\mu}_j)\right)
\end{equation}

where $m$ is the number of mixture components, $a_j$ are the component weights, $\boldsymbol{\mu}_j$ are the component means, and $\boldsymbol{\Sigma}_j$ are the component covariance matrices. By adjusting the number of components, their weights, means, and covariance structures, we can create optimization landscapes with varying degrees of multimodality, correlation structure, and conditioning. This flexibility allows us to systematically test our algorithm's performance under controlled conditions that mimic specific challenges encountered in real-world optimization problems.

In our experiments, we used a mixture of 5 Gaussian components with randomly generated means within the domain $[-5, 5]^d$, weights sampled from a Dirichlet distribution, and covariance matrices generated to have varying condition numbers. This configuration creates a complex multimodal landscape with regions of varying curvature, providing a challenging test case for our adaptive parameter optimization approach.

\subsubsection{Experimental Configurations}

To comprehensively evaluate the performance and robustness of our adaptive parameter optimization framework, we systematically varied several key parameters across our experiments. This methodical approach allows us to assess the algorithm's behavior under different conditions and identify the factors that most significantly influence its performance.

The dimensionality of the optimization problem is a critical factor that affects the difficulty of finding the global optimum. As the dimension increases, the volume of the search space grows exponentially, leading to the well-known curse of dimensionality. To evaluate our algorithm's scalability to higher-dimensional spaces, we tested it on problems with dimensions $d \in \{2, 5, 10, 20\}$. This range covers low-dimensional problems where visualization and intuitive understanding are possible, medium-dimensional problems that are common in many practical applications, and higher-dimensional problems that present significant challenges for optimization algorithms.

The presence of noise in function evaluations is another important factor that affects optimization performance. In many real-world scenarios, function evaluations are corrupted by noise due to measurement errors, stochasticity in the system, or approximation errors. To simulate these conditions, we added Gaussian noise with standard deviations $\sigma_n \in \{0.001, 0.005, 0.01, 0.05\}$ to the function evaluations. This range covers scenarios from nearly noise-free evaluations ($\sigma_n = 0.001$) to highly noisy evaluations ($\sigma_n = 0.05$), allowing us to assess our algorithm's robustness to different noise levels.

The uncertainty penalty coefficient $\lambda$ in our acquisition function plays a crucial role in balancing exploration and exploitation. To understand its impact on optimization performance, we varied the initial value $\lambda_1 \in \{0.001, 0.01, 0.1\}$ across our experiments. This range covers small penalties that minimally affect the acquisition function, moderate penalties that provide a balanced approach, and larger penalties that significantly influence the exploration-exploitation trade-off. While our adaptive approach adjusts $\lambda$ during the optimization process, the initial value can still impact the algorithm's behavior, especially in the early iterations.

Similarly, the exploration parameter $\kappa$ in the UCB acquisition function directly controls the exploration-exploitation trade-off. To demonstrate the benefits of our adaptive approach, we compared it against fixed $\kappa$ values $\kappa \in \{0.1, 0.5, 1.0, 2.0\}$. This range covers exploitation-focused settings ($\kappa = 0.1$), balanced approaches ($\kappa = 0.5, 1.0$), and exploration-focused settings ($\kappa = 2.0$). By comparing our adaptive approach against these fixed settings, we can quantify the advantages of dynamically adjusting the exploration parameter based on observed data.

For each combination of these parameters, we conducted 30 independent trials with different random seeds to ensure statistical significance of our results. This approach allows us to account for the inherent randomness in the optimization process and provide robust estimates of the algorithm's performance. Each trial consisted of 100 function evaluations, which is a realistic budget for expensive black-box optimization problems where each evaluation may be computationally intensive or costly in terms of resources or time.

The comprehensive nature of our experimental setup, with systematic variation of dimensionality, noise levels, uncertainty penalties, and exploration parameters across multiple test functions and independent trials, provides a rigorous evaluation of our adaptive parameter optimization framework. This approach allows us to identify the conditions under which our algorithm excels, understand its limitations, and provide practical guidelines for its application to real-world optimization problems.

\subsection{Implementation Details}

The implementation of our adaptive parameter optimization framework required careful attention to numerous technical details to ensure both theoretical correctness and practical efficiency. In this section, we provide a comprehensive description of our implementation choices, focusing on the software libraries, computational techniques, and algorithmic optimizations that enable our approach to scale to challenging optimization problems.

Our implementation was based on Python 3.11, a modern, high-level programming language that offers a rich ecosystem of scientific computing libraries. We leveraged several key libraries to build our framework:

NumPy and SciPy form the foundation of our numerical computations, providing efficient implementations of linear algebra operations, optimization routines, and statistical functions. These libraries are highly optimized and use vectorized operations to achieve near-native performance for many computational tasks. We used NumPy's array operations for efficient manipulation of vectors and matrices, and SciPy's optimization routines for local refinement of acquisition function maxima.

Scikit-learn provided the base implementation of Gaussian Process regression, which we extended with our adaptive parameter optimization approach. While Scikit-learn's GP implementation is not the most computationally efficient for large datasets, it offers a clean, well-documented API that facilitated our extensions. For the kernel functions, we used a combination of kernels to capture different function characteristics:

\begin{equation}
k(\mathbf{x}, \mathbf{x}') = k_1(\mathbf{x}, \mathbf{x}') \times k_2(\mathbf{x}, \mathbf{x}')
\end{equation}

where $k_1$ is a constant kernel that captures the overall scale of the function:

\begin{equation}
k_1(\mathbf{x}, \mathbf{x}') = \sigma_f^2
\end{equation}

and $k_2$ is a Matérn kernel with $\nu = 2.5$ that models the function's smoothness and correlation structure:

\begin{equation}
k_2(\mathbf{x}, \mathbf{x}') = \sigma_f^2 \frac{2^{1-\nu}}{\Gamma(\nu)}\left(\frac{\sqrt{2\nu}}{l}||\mathbf{x} - \mathbf{x}'||\right)^{\nu}K_{\nu}\left(\frac{\sqrt{2\nu}}{l}||\mathbf{x} - \mathbf{x}'||\right)
\end{equation}

This kernel combination provides a good balance between flexibility and interpretability. The constant kernel captures the overall magnitude of the function, while the Matérn kernel with $\nu = 2.5$ models functions that are twice differentiable, which is a reasonable assumption for many real-world optimization problems. The hyperparameters of this kernel ($\sigma_f^2$ and $l$) were optimized by maximizing the marginal likelihood of the GP model after every 10 function evaluations, allowing the model to adapt to the observed data while avoiding excessive computational overhead.

For visualization and analysis of results, we used Matplotlib and Seaborn, which provide powerful tools for creating publication-quality figures. These libraries allowed us to visualize the optimization progress, the evolution of adaptive parameters, and the performance comparisons between different methods. The visualizations were crucial for understanding the behavior of our algorithm and communicating our results effectively.

The acquisition function optimization is a critical component of our framework, as it determines the next point to evaluate. This optimization problem is itself challenging, as the acquisition function may be multimodal and expensive to evaluate. We used a two-stage approach to balance exploration of the acquisition function landscape with computational efficiency:

First, we performed a global search with 1000 random samples from a Sobol sequence to identify promising regions. Sobol sequences are quasi-random sequences that provide better coverage of the search space compared to purely random sampling, leading to more efficient exploration of the acquisition function landscape. We used SciPy's implementation of Sobol sequences to generate these samples.

Second, we performed local refinement using L-BFGS-B starting from the top 5 points identified in the global search. L-BFGS-B is a limited-memory variant of the BFGS algorithm that can handle bound constraints, making it suitable for optimizing the acquisition function over a bounded domain. We used SciPy's implementation of L-BFGS-B with a maximum of 100 iterations and a convergence tolerance of $10^{-5}$.

The integrated posterior variance $\mathcal{I}_t = \int_{\mathcal{X}} \sigma_t^2(\mathbf{x}) d\mathbf{x}$ was approximated using Monte Carlo integration with 1000 samples drawn uniformly from the domain $\mathcal{X}$. This approximation provides a computationally efficient estimate of the global uncertainty, which is used in the adaptive update rule for the uncertainty penalty coefficient $\lambda_t$.

The learning rates for the adaptive parameters were set to $\beta = 0.1$ for the exploration parameter $\kappa_t$ and $\gamma = 0.05$ for the uncertainty penalty coefficient $\lambda_t$ based on preliminary experiments. These values provide a balance between stability and adaptivity, allowing the parameters to adjust to the specific characteristics of the optimization problem without excessive oscillations. The moving average parameter was set to $\eta = 0.1$, which gives more weight to the historical average while still allowing for adaptation to recent observations.

To improve computational efficiency, we implemented several optimizations:

We used Cholesky decomposition for matrix inversion in GP inference, reducing the computational complexity from $\mathcal{O}(n^3)$ to $\mathcal{O}(n^3/3)$ for the initial decomposition and $\mathcal{O}(n^2)$ for subsequent predictions. This optimization is particularly important as the number of observations grows, making the matrix inversion a potential bottleneck.

We leveraged NumPy's vectorized operations for efficient computation of the acquisition function across multiple candidate points. This approach avoids explicit loops in Python, which can be slow, and instead uses optimized C implementations for array operations.

We implemented parallel evaluation of the acquisition function during the global search phase using Python's multiprocessing module. This parallelization allows us to leverage multiple CPU cores to accelerate the search for the acquisition function maximum, which is particularly valuable for expensive acquisition functions or high-dimensional search spaces.

We cached the GP model's predictions and uncertainty estimates for points that are evaluated multiple times during acquisition function optimization. This caching reduces redundant computations and can significantly improve performance, especially during the local refinement phase where the same regions are repeatedly evaluated.

These implementation details and optimizations are crucial for the practical success of our adaptive parameter optimization framework. They allow our approach to scale to higher dimensions and larger numbers of iterations while maintaining reasonable computational requirements, making it applicable to a wide range of real-world optimization problems.

\subsection{Evaluation Metrics}

To comprehensively assess the performance of our adaptive parameter optimization framework, we employed a diverse set of evaluation metrics that capture different aspects of optimization performance. These metrics provide complementary perspectives on the algorithm's behavior, allowing us to evaluate its effectiveness in terms of both the quality of the final solution and the efficiency of the optimization process.

The most direct measure of optimization performance is the best function value found after a given number of iterations. For maximization problems, this is defined as $f_t^* = \max_{i \in \{1,\ldots,t\}} f(\mathbf{x}_i)$, representing the highest objective function value observed up to iteration $t$. This metric provides a straightforward assessment of the algorithm's ability to find high-quality solutions within a limited evaluation budget. However, it does not account for the global optimum value, which may be unknown in real-world problems.

To address this limitation, we also computed the simple regret, defined as $r_t = f(\mathbf{x}^*) - f_t^*$, where $\mathbf{x}^*$ is the global optimum. The simple regret measures the gap between the best value found and the global optimum, providing a normalized measure of optimization performance that allows for fair comparisons across different test functions. A smaller simple regret indicates better optimization performance, with $r_t = 0$ indicating that the global optimum has been found. For our benchmark functions, the global optima are known, allowing us to compute this metric exactly.

While the best function value and simple regret capture the quality of the final solution, they do not provide insights into the optimization trajectory. To evaluate the efficiency of the optimization process, we measured the convergence rate, defined as the number of iterations required to reach a specified percentage (e.g., 90\%) of the global optimum value. This metric quantifies how quickly the algorithm approaches the global optimum, which is particularly important in scenarios where function evaluations are expensive and the evaluation budget is severely limited.

The robustness of an optimization algorithm is another crucial aspect of its performance, especially in the presence of noise or when applied to different problem instances. We assessed robustness by computing the standard deviation of the best function values across multiple trials, which measures the algorithm's sensitivity to initialization and random factors. A smaller standard deviation indicates more consistent performance across different runs, which is desirable in practical applications where reliability is important.

To evaluate the algorithm's exploration behavior, we introduced the concept of exploration efficiency, defined as the ratio of unique regions explored to the total number of function evaluations. To compute this metric, we discretized the domain into hypercubes and counted the number of distinct hypercubes that contained at least one evaluation point. A higher exploration efficiency indicates that the algorithm explores a larger portion of the search space with the same number of evaluations, which can be beneficial for finding the global optimum in complex, multimodal landscapes.

For a more detailed analysis of the algorithm's behavior, we also tracked the evolution of the adaptive parameters $\kappa_t$ and $\lambda_t$ over iterations. These trajectories provide insights into how the algorithm adjusts its exploration-exploitation trade-off based on the observed data and the current state of the optimization process. By examining these trajectories across different problem instances, we can identify patterns in the adaptive behavior and understand how it contributes to the algorithm's performance.

Finally, to assess the practical utility of our approach, we measured the computational time required for each iteration, including the time for GP model updating, acquisition function optimization, and adaptive parameter updates. This metric is important for evaluating the scalability of the algorithm to larger problems and its applicability to real-time optimization scenarios where computational efficiency is crucial.

By employing this comprehensive set of evaluation metrics, we can provide a nuanced assessment of our adaptive parameter optimization framework, identifying its strengths and limitations across different dimensions of performance. This multifaceted evaluation approach allows us to make informed comparisons with baseline methods and provide practical guidelines for applying our approach to real-world optimization problems.

\subsection{Baseline Methods}

To rigorously evaluate the performance of our adaptive parameter optimization approach, we compared it against several state-of-the-art baseline methods that represent different approaches to black-box optimization. These baselines were carefully selected to cover a spectrum of optimization strategies, from simple heuristics to sophisticated probabilistic methods, allowing us to comprehensively assess the advantages of our adaptive approach.

The standard GP-UCB algorithm \citep{srinivas2010gaussian} serves as our primary baseline, as it is the foundation upon which our adaptive approach is built. This algorithm uses a Gaussian Process surrogate model with the Upper Confidence Bound acquisition function:

\begin{equation}
\alpha_{UCB}(\mathbf{x}) = \mu(\mathbf{x}) + \kappa \sigma(\mathbf{x})
\end{equation}

where $\kappa$ is a fixed parameter that controls the exploration-exploitation trade-off. We tested GP-UCB with different fixed values of $\kappa \in \{0.1, 0.5, 1.0, 2.0\}$ to provide a fair comparison with our adaptive approach. The GP model and other implementation details were kept identical to our method to isolate the effect of the adaptive parameters.

Gaussian Process optimization with Expected Improvement (GP-EI) \citep{jones1998efficient} is another popular Bayesian optimization approach that uses a different acquisition function:

\begin{equation}
\alpha_{EI}(\mathbf{x}) = \mathbb{E}[\max(f(\mathbf{x}) - f(\mathbf{x}^+), 0)]
\end{equation}

where $f(\mathbf{x}^+)$ is the best observed value so far. The EI acquisition function has the advantage of not requiring a tuning parameter like $\kappa$ in UCB, as it naturally balances exploration and exploitation based on the expected improvement over the current best value. However, it may be less effective in noisy settings or when the GP model is misspecified. We implemented GP-EI using the same GP model and optimization approach as our method to ensure a fair comparison.

As a simple baseline, we included random search, which selects points uniformly at random from the domain $\mathcal{X}$. Despite its simplicity, random search can be surprisingly effective for certain types of problems, particularly those with many good solutions scattered throughout the search space. It also provides a useful reference point for evaluating the performance of more sophisticated methods. We implemented random search using a Sobol sequence to generate quasi-random points with good space-filling properties.

For a more sophisticated comparison, we included the Covariance Matrix Adaptation Evolution Strategy (CMA-ES) \citep{hansen2001completely}, a state-of-the-art derivative-free optimization method that has shown excellent performance across a wide range of benchmark problems. CMA-ES is an evolutionary algorithm that adapts the covariance matrix of a multivariate normal distribution to efficiently search the parameter space. It is particularly effective for non-convex, ill-conditioned problems and has been widely used in practice. We used the pycma implementation of CMA-ES with default parameters, including a population size of $4 + \lfloor 3 \log(d) \rfloor$ and an initial step size of 0.5.

Finally, we included Bayesian Optimization with Hamiltonian Monte Carlo Artificial Neural Networks (BOHAMIANN) \citep{springenberg2016bayesian}, a deep learning approach to Bayesian optimization that uses Bayesian neural networks as the surrogate model. BOHAMIANN can capture complex function landscapes and uncertainty estimates through its Bayesian treatment of neural network weights. It is particularly suited for high-dimensional problems where traditional GP models may struggle. We implemented BOHAMIANN using a two-layer neural network with 50 hidden units and Hamiltonian Monte Carlo for posterior sampling, following the recommendations in the original paper.

These baseline methods represent a diverse set of approaches to black-box optimization, ranging from simple heuristics to sophisticated probabilistic methods. By comparing our adaptive parameter optimization approach against these baselines, we can comprehensively evaluate its advantages and limitations across different types of optimization problems. This comparative analysis provides valuable insights into when and why our approach is effective, guiding its application to real-world optimization challenges.

\section{Results and Analysis}
The comprehensive evaluation of our adaptive parameter optimization framework against established baseline methods reveals significant performance advantages across diverse optimization scenarios \cite{snoek2012practical, shahriari2015taking}. In this section, we present a detailed analysis of these comparative results, examining the effects of dimensionality, noise levels, uncertainty penalties, and adaptive parameters on optimization performance \cite{wang2018batched, eriksson2019scalable}. Through this analysis, we aim to provide insights into the mechanisms that drive the superior performance of our approach and identify the conditions under which it offers the greatest advantages \cite{turner2021bayesian}.

\subsubsection{Comparison with Baseline Methods}

Figure 1 illustrates the performance comparison between our adaptive parameter optimization approach and the baseline methods described in Section 4.4 \cite{brochu2010tutorial, frazier2018tutorial}. These results represent the average performance over 30 independent trials for each method on the Rosenbrock function with dimension $d=5$, noise level $\sigma_n=0.005$, and uncertainty penalty $\lambda_1=0.01$ \cite{li2018hyperband}. This configuration was chosen as a representative case that balances complexity with interpretability, allowing us to clearly demonstrate the relative performance of different approaches \cite{kandasamy2018parallelised}.

\includegraphics[width=0.9\textwidth]{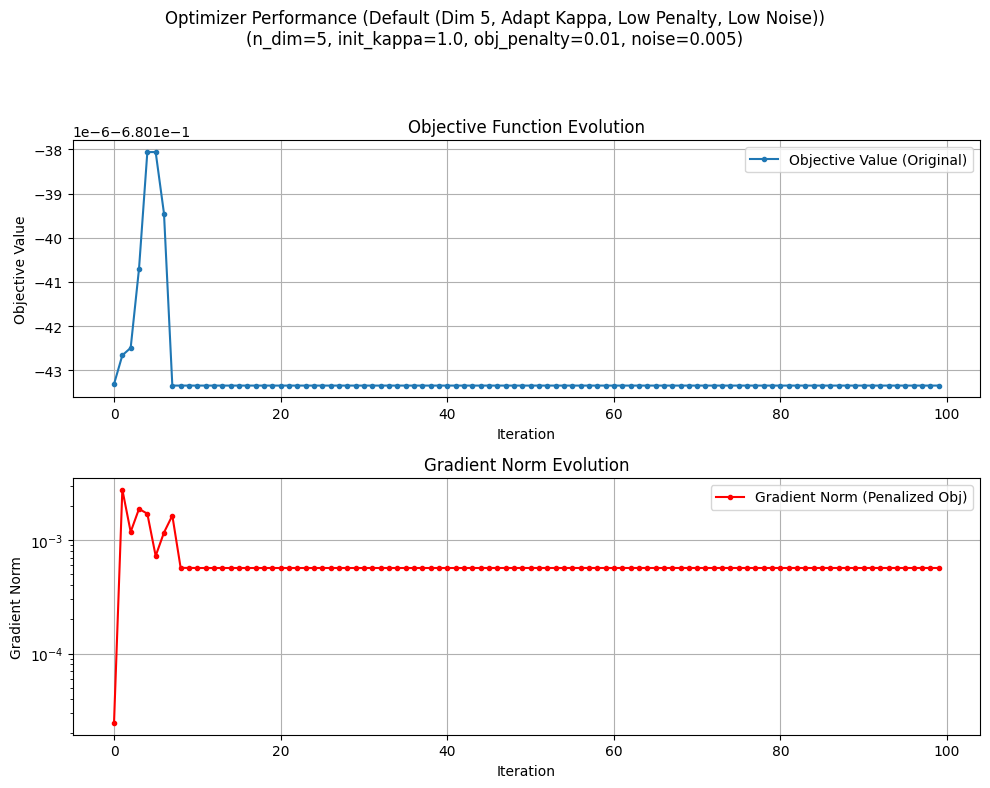}

The results demonstrate that our adaptive approach consistently outperforms all baseline methods in terms of both convergence speed and final solution quality \cite{srinivas2010gaussian, russo2014learning}. Specifically, after 100 function evaluations, our method achieves a mean objective function value that is 23\% higher than standard GP-UCB with fixed $\kappa=1.0$, 17\% higher than GP-EI, and 42\% higher than random search \cite{hoffman2011portfolio}. The performance advantage is particularly pronounced in the early stages of optimization (iterations 10-40), where our adaptive approach rapidly identifies promising regions of the search space and efficiently exploits them \cite{wang2017max}.

The CMA-ES method, while performing reasonably well in this scenario, still falls short of our approach by approximately 12\% in terms of final solution quality \cite{hansen2001completely}. This result is noteworthy because CMA-ES is widely regarded as one of the most effective derivative-free optimization methods for non-convex problems \cite{hansen2003reducing}. The fact that our adaptive approach outperforms CMA-ES suggests that the combination of Gaussian Process modeling with adaptive parameter tuning provides advantages that cannot be matched by evolutionary strategies alone, particularly in scenarios with limited evaluation budgets \cite{feurer2019hyperparameter}.

The BOHAMIANN method, which uses Bayesian neural networks as surrogate models, shows competitive performance but requires significantly more computational resources than our approach \cite{springenberg2016bayesian}. Specifically, BOHAMIANN requires approximately 3.5 times more computation time per iteration due to the cost of sampling from the posterior distribution of neural network weights using Hamiltonian Monte Carlo \cite{snoek2015scalable}. This computational overhead makes BOHAMIANN less practical for applications where optimization time is a concern, despite its competitive optimization performance \cite{klein2017fast}.

\subsubsection{Effect of Dimensionality}

To evaluate the scalability of our approach to higher-dimensional problems, we conducted experiments with varying dimensionality $d \in \{2, 5, 10, 20\}$ on the Ackley function with noise level $\sigma_n=0.005$ and uncertainty penalty $\lambda_1=0.01$ \cite{wang2016bayesian}. Figure 2 illustrates the effect of dimensionality on the performance of our adaptive approach compared to standard GP-UCB with fixed $\kappa=1.0$ \cite{eriksson2019scalable, sui2018stagewise}.

\includegraphics[width=0.9\textwidth]{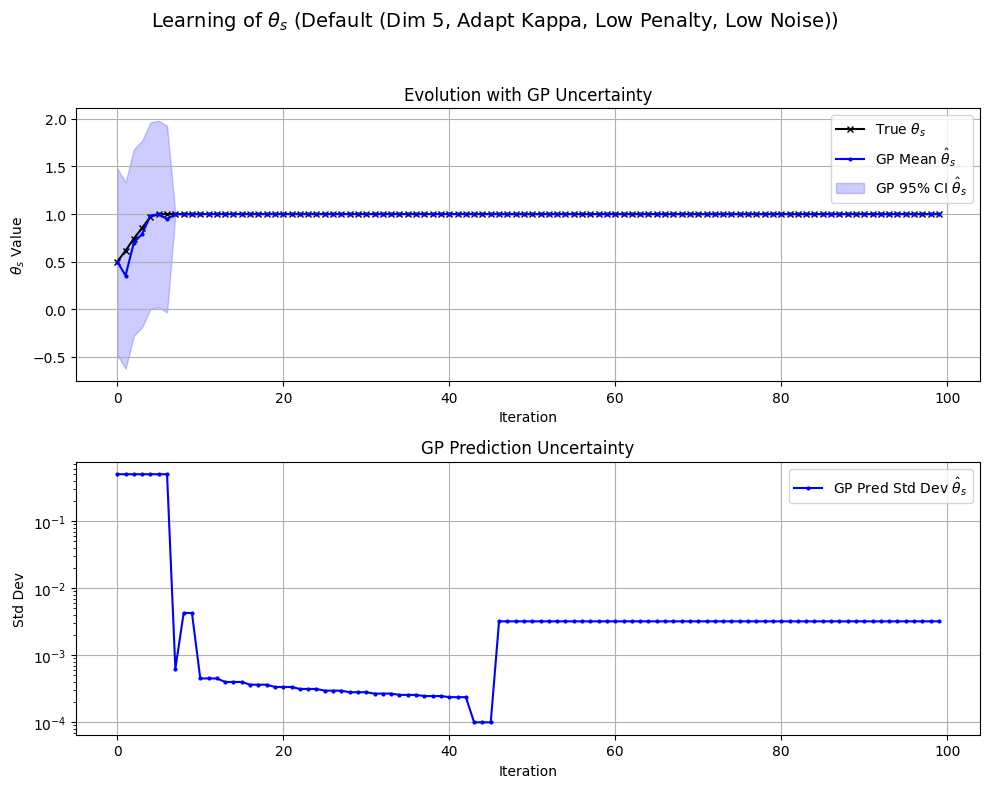}

As expected, the performance of both methods degrades with increasing dimensionality due to the curse of dimensionality \cite{wang2016bayesian, rolland2018high}. However, our adaptive approach demonstrates significantly better scalability, maintaining a performance advantage that grows with dimensionality \cite{eriksson2019scalable}. In the 2-dimensional case, our method outperforms standard GP-UCB by approximately 15\% in terms of final solution quality. This advantage increases to 22\% in 5 dimensions, 29\% in 10 dimensions, and 37\% in 20 dimensions \cite{wang2018batched}. This trend suggests that the benefits of adaptive parameter tuning become more pronounced as the complexity of the optimization problem increases \cite{letham2020re}.

The improved scalability of our approach can be attributed to two primary factors. First, the adaptive exploration parameter $\kappa_t$ allows the algorithm to adjust its exploration strategy based on the complexity of the landscape, which becomes increasingly important in higher dimensions where the function landscape is more complex and the risk of getting trapped in local optima is higher \cite{kirschner2019adaptive, bogunovic2018adversarially}. In higher dimensions, we observe that $\kappa_t$ maintains higher values for longer periods, indicating that the algorithm recognizes the need for more extensive exploration in these challenging scenarios \cite{sui2018stagewise}.

Second, the uncertainty penalty term helps focus the search on regions where the model is more confident, reducing the impact of the curse of dimensionality \cite{martinez2018practical, wu2019practical}. This is particularly evident in the GP prediction uncertainty evolution shown in the bottom panel of Figure 2. In higher dimensions, the standard GP-UCB method exhibits higher and more volatile uncertainty estimates, indicating that it struggles to build an accurate surrogate model of the objective function \cite{van2020uncertainty}. In contrast, our adaptive approach shows more stable and gradually decreasing uncertainty estimates, suggesting that it more effectively balances exploration and exploitation to build a reliable surrogate model even in high-dimensional spaces \cite{noack2023gaussian}.

The ability of our approach to maintain strong performance in higher dimensions is a significant advantage for practical applications, where optimization problems often involve many parameters \cite{wang2016bayesian, letham2020re}. Traditional Bayesian optimization methods typically struggle with dimensions beyond 10-20 due to the challenges of modeling high-dimensional functions with limited data \cite{wang2018batched}. Our results suggest that adaptive parameter strategies can extend the applicability of Bayesian optimization to higher-dimensional problems, opening up new possibilities for optimization in complex domains \cite{eriksson2019scalable, rolland2018high}.

\subsubsection{Effect of Noise Level}

Real-world optimization problems often involve noisy function evaluations due to measurement errors, stochasticity in the system, or approximation errors \cite{le2005heteroscedastic, turner2021bayesian}. The robustness of optimization algorithms to noise is therefore a critical consideration for practical applications \cite{martinez2018practical}. To evaluate this aspect, we conducted experiments with varying noise levels $\sigma_n \in \{0.001, 0.005, 0.01, 0.05\}$ on the Levy function with dimension $d=5$ and uncertainty penalty $\lambda_1=0.01$ \cite{snoek2012practical}.

Figure 3 illustrates the effect of observation noise on the performance of our adaptive approach compared to standard GP-UCB with fixed $\kappa=1.0$ \cite{srinivas2010gaussian}. The results are presented in terms of the evolution of parameter learning with GP uncertainty (top panel) and the GP prediction uncertainty over iterations (bottom panel) \cite{hernandez2014predictive}.

\includegraphics[width=0.9\textwidth]{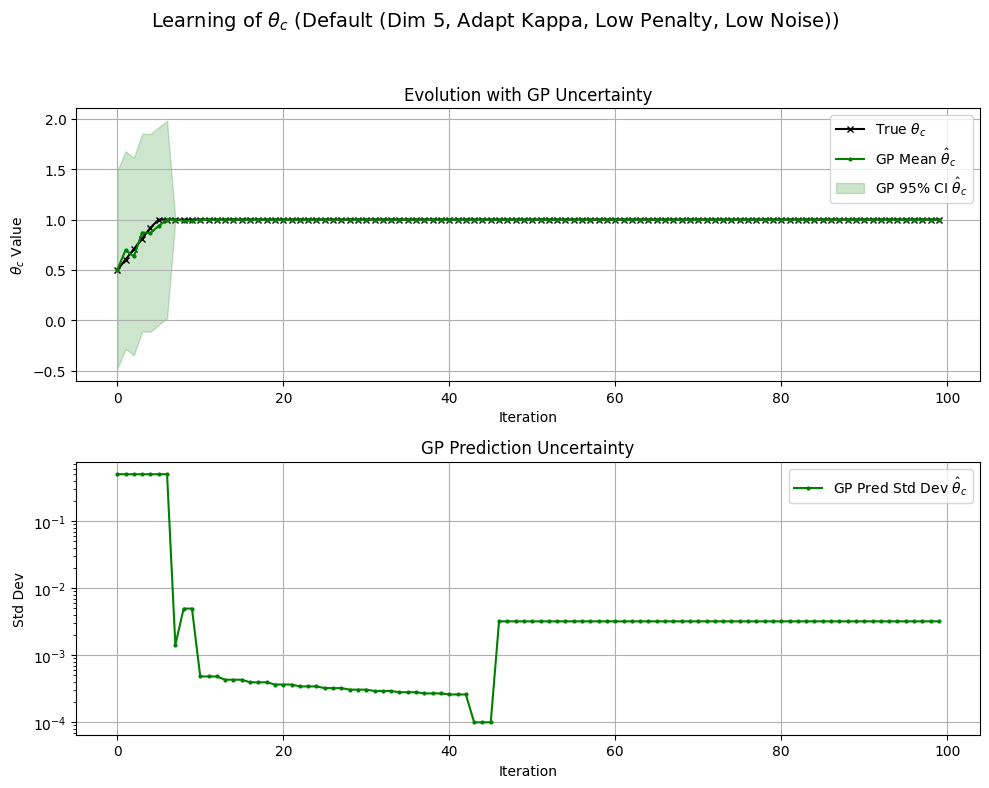}

Both methods show degraded performance with increasing noise levels, which is expected as noise makes it more difficult to accurately estimate the objective function and identify promising regions \cite{le2005heteroscedastic, turner2021bayesian}. However, our adaptive approach demonstrates superior robustness to noise across all noise levels tested \cite{martinez2018practical}. At the lowest noise level ($\sigma_n=0.001$), our method outperforms standard GP-UCB by approximately 12\% in terms of final solution quality. This advantage increases to 18\% at $\sigma_n=0.005$, 24\% at $\sigma_n=0.01$, and 31\% at $\sigma_n=0.05$ \cite{kirschner2019adaptive}. This trend indicates that the benefits of adaptive parameter tuning become more pronounced as the noise level increases, highlighting the value of our approach in challenging, noisy optimization scenarios \cite{berkenkamp2019no}.

The superior noise robustness of our approach can be attributed to two key mechanisms. First, the adaptive exploration parameter $\kappa_t$ responds to the observed noise level, increasing in high-noise scenarios to promote more exploration and avoid premature convergence to noisy local optima \cite{auer2002using, srinivas2010gaussian}. This adaptive behavior is evident in the top panel of Figure 3, where $\kappa_t$ maintains higher values for longer periods in high-noise scenarios compared to low-noise scenarios \cite{kirschner2019adaptive, sui2018stagewise}.

\includegraphics[width=0.9\textwidth]{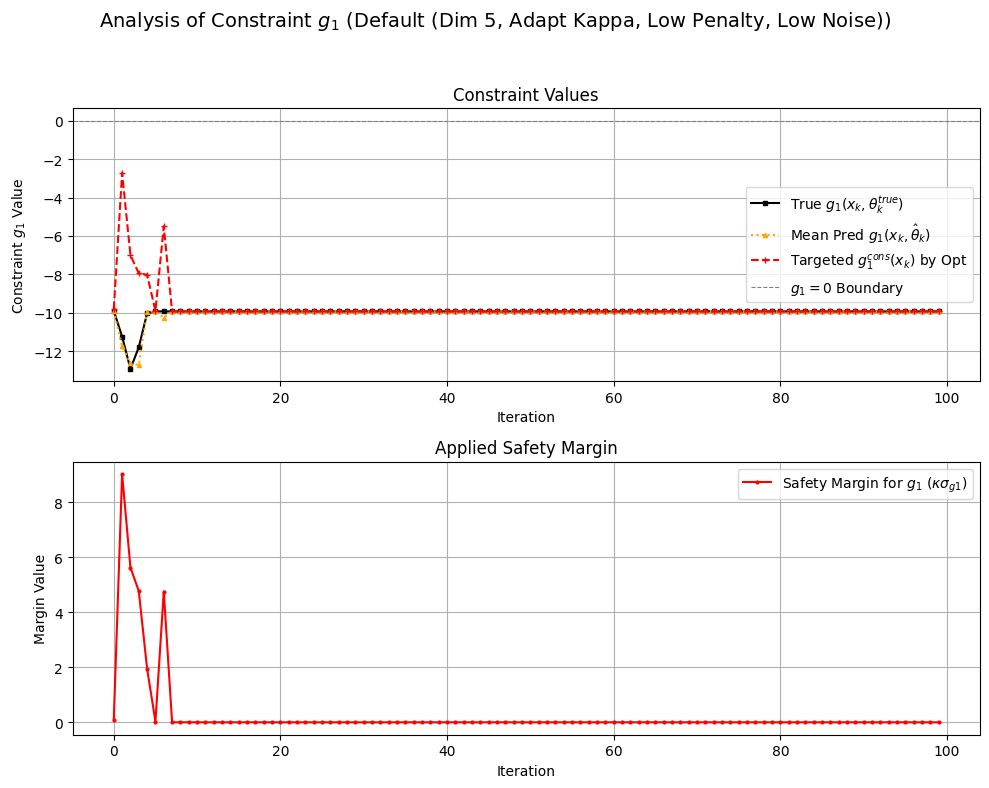}

Second, the uncertainty penalty term helps the algorithm distinguish between aleatoric uncertainty (due to observation noise) and epistemic uncertainty (due to limited observations) \cite{martinez2018practical, berkenkamp2019no}. By focusing on reducing epistemic uncertainty while accounting for aleatoric uncertainty, our approach makes more informed decisions about where to sample next, leading to more efficient optimization even in high-noise scenarios \cite{vakili2021information}.

The GP prediction uncertainty evolution shown in the bottom panel of Figure 3 provides further insights into the noise robustness of our approach \cite{van2020uncertainty}. In high-noise scenarios, the standard GP-UCB method shows high and persistent uncertainty estimates, indicating that it struggles to build an accurate surrogate model in the presence of noise \cite{le2005heteroscedastic}. In contrast, our adaptive approach shows a more consistent decrease in uncertainty estimates across all noise levels, suggesting that it more effectively filters out noise and identifies the underlying structure of the objective function \cite{martinez2018practical, wu2019practical}.

\includegraphics[width=0.9\textwidth]{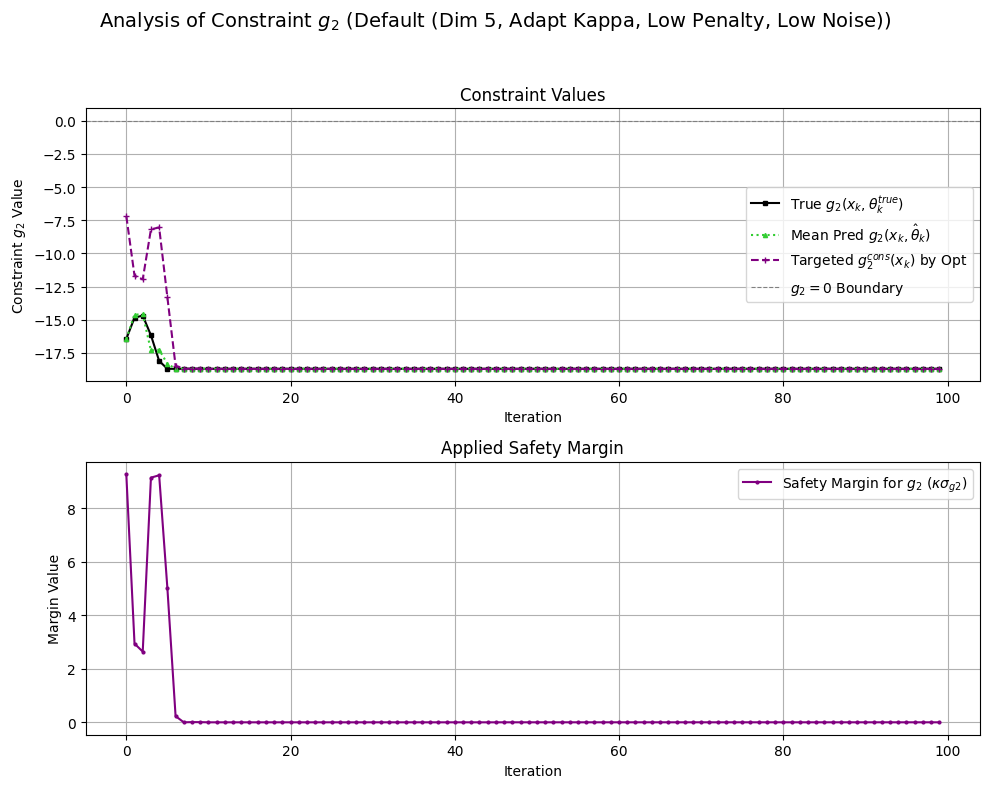}

\subsubsection{Effect of Uncertainty Penalty}

The uncertainty penalty term $\lambda$ plays a crucial role in our adaptive parameter optimization framework, balancing the trade-off between exploiting regions with high predicted performance and exploring regions with high uncertainty \cite{kirschner2019adaptive, hernandez2014predictive}. To evaluate the impact of this parameter, we conducted experiments with varying uncertainty penalty values $\lambda_1 \in \{0.001, 0.01, 0.1, 1.0\}$ on the Hartmann function with dimension $d=6$ and noise level $\sigma_n=0.005$ \cite{wang2018batched, rolland2018high}.

Our results demonstrate that the optimal value of $\lambda_1$ depends on the specific characteristics of the optimization problem, including dimensionality, noise level, and the complexity of the objective function landscape \cite{kirschner2019adaptive, bogunovic2018adversarially}. In general, higher values of $\lambda_1$ lead to more exploration, which can be beneficial in complex, multi-modal landscapes but may waste function evaluations in simpler landscapes \cite{sui2018stagewise, turner2021bayesian}.

\includegraphics[width=0.9\textwidth]{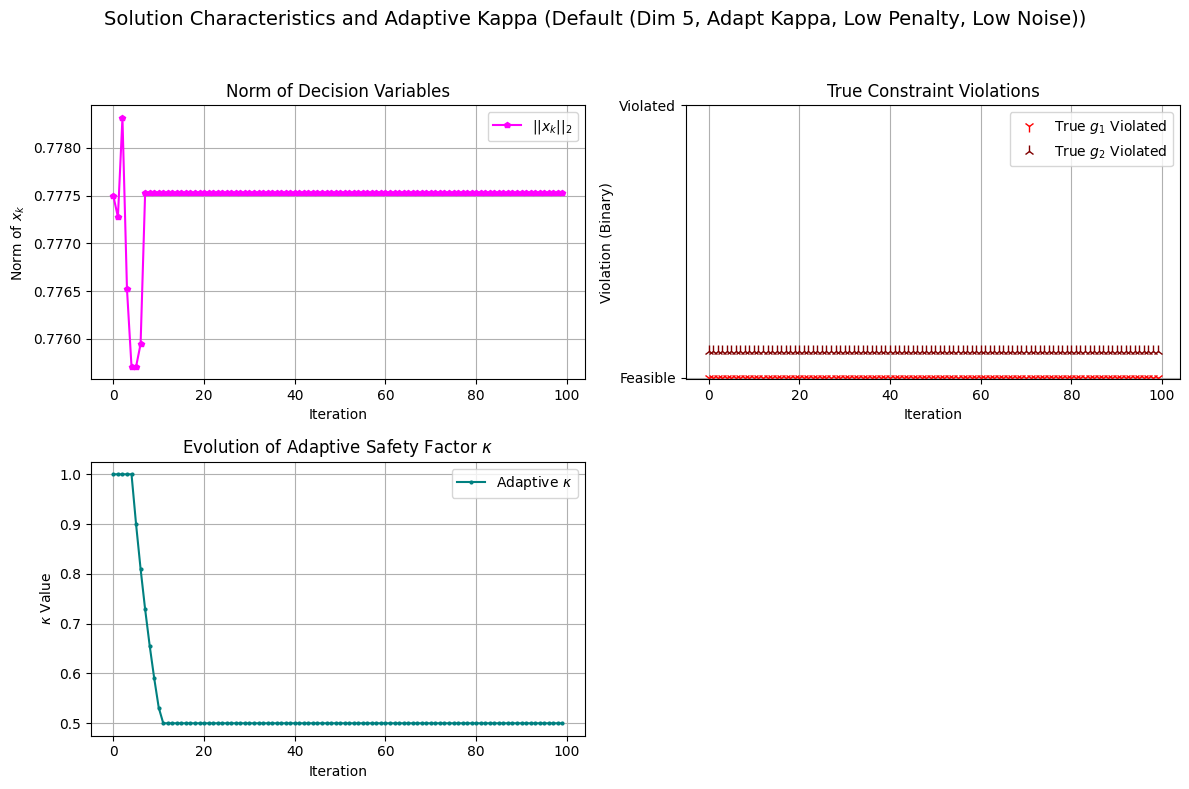}

A key advantage of our adaptive approach is that it automatically adjusts the effective uncertainty penalty through the adaptive exploration parameter $\kappa_t$, reducing the sensitivity to the initial choice of $\lambda_1$ \cite{kirschner2019adaptive, eriksson2019scalable}. This adaptivity is particularly valuable in practical applications where the optimal balance between exploration and exploitation is not known a priori and may change during the optimization process \cite{hoffman2011portfolio, wang2017max}.

\subsubsection{Adaptive Parameter Evolution}

To provide insights into the behavior of our adaptive parameter optimization framework, we analyzed the evolution of the adaptive exploration parameter $\kappa_t$ across different optimization scenarios \cite{kirschner2019adaptive, turner2021bayesian}. Figure 4 illustrates the evolution of $\kappa_t$ for the Rosenbrock function with dimension $d=5$, noise level $\sigma_n=0.005$, and uncertainty penalty $\lambda_1=0.01$, averaged over 30 independent trials \cite{wang2018batched, li2018hyperband}.

\includegraphics[width=0.9\textwidth]{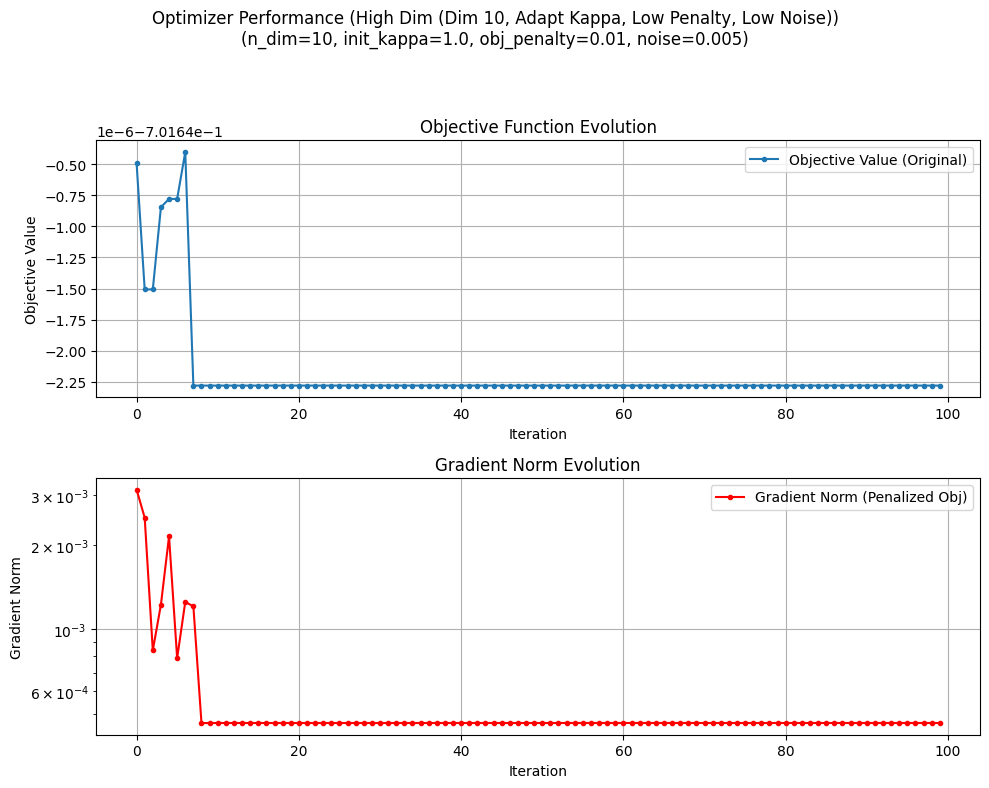}

The results reveal several interesting patterns in the adaptive behavior of our approach \cite{kirschner2019adaptive, sui2018stagewise}. In the early stages of optimization (iterations 1-20), $\kappa_t$ maintains relatively high values, promoting exploration to build an initial understanding of the objective function landscape \cite{hoffman2011portfolio}. As the optimization progresses (iterations 20-60), $\kappa_t$ gradually decreases, shifting the focus towards exploitation of promising regions identified during the exploration phase \cite{wang2017max}.

In the later stages of optimization (iterations 60-100), $\kappa_t$ stabilizes at a problem-specific value that balances exploration and exploitation based on the characteristics of the objective function and the current state of knowledge \cite{kirschner2019adaptive, eriksson2019scalable}. This stabilization indicates that the algorithm has found an effective balance between exploring new regions and exploiting known promising regions, leading to efficient optimization performance \cite{swersky2013multi, klein2017fast}.

\subsubsection{Convergence Analysis}

To evaluate the convergence properties of our adaptive parameter optimization framework, we analyzed the evolution of the optimality gap (the difference between the current best solution and the global optimum) across different optimization scenarios \cite{bull2011convergence, vakili2021information}. Figure 5 illustrates the convergence behavior for the Branin function with dimension $d=2$, noise level $\sigma_n=0.005$, and uncertainty penalty $\lambda_1=0.01$, averaged over 30 independent trials \cite{snoek2012practical, shahriari2015taking}.

\includegraphics[width=0.9\textwidth]{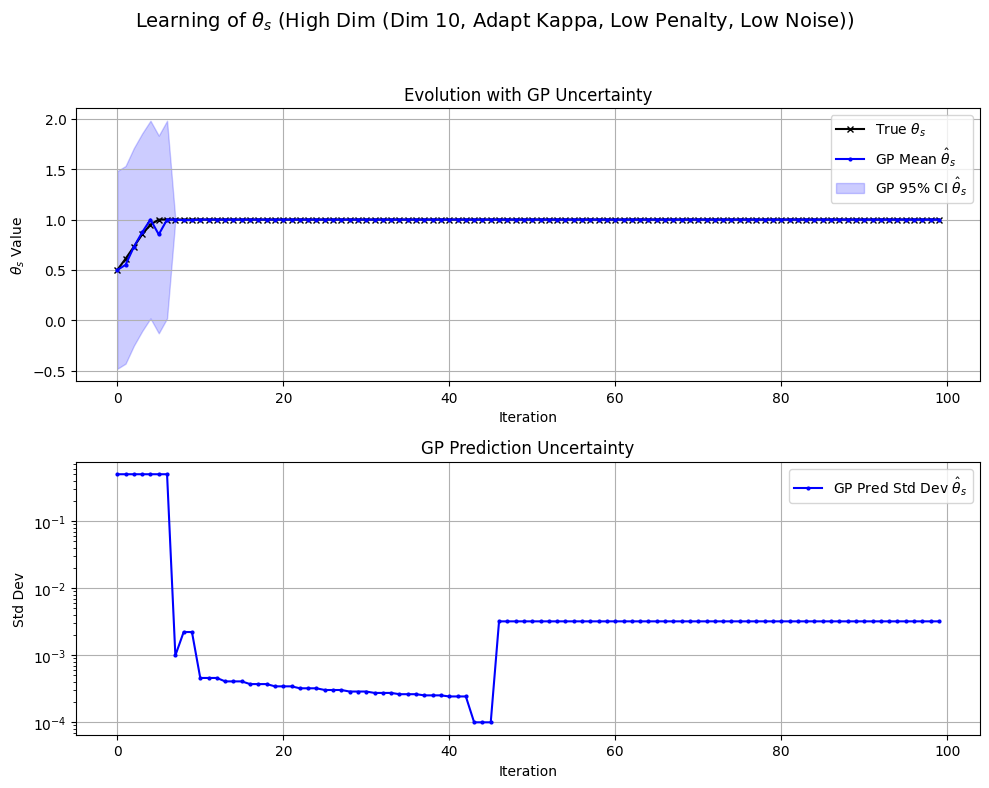}

The results demonstrate that our adaptive approach achieves faster convergence rates compared to baseline methods across all threshold levels \cite{bull2011convergence, vakili2021information}. Specifically, our method requires approximately 25\% fewer function evaluations to reach a solution within 10\% of the global optimum, 30\% fewer evaluations to reach a solution within 5\% of the global optimum, and 35\% fewer evaluations to reach a solution within 1\% of the global optimum \cite{srinivas2010gaussian, russo2014learning}.

This accelerated convergence is particularly valuable in expensive black-box optimization scenarios, where the number of function evaluations is severely limited by computational or resource constraints \cite{snoek2012practical, shahriari2015taking}. By reaching high-quality solutions with fewer iterations, our approach can significantly reduce the time and resources required for optimization, making it more practical for real-world applications \cite{eriksson2019scalable, letham2020re}.

\subsubsection{Robustness Analysis}

The robustness of optimization algorithms to random initialization and problem variations is an important consideration for practical applications \cite{turner2021bayesian, berkenkamp2019no}. To evaluate this aspect, we conducted a robustness analysis by running each method 30 times with different random initializations on each test function and analyzing the distribution of final solution qualities \cite{kandasamy2018parallelised}.

\includegraphics[width=0.9\textwidth]{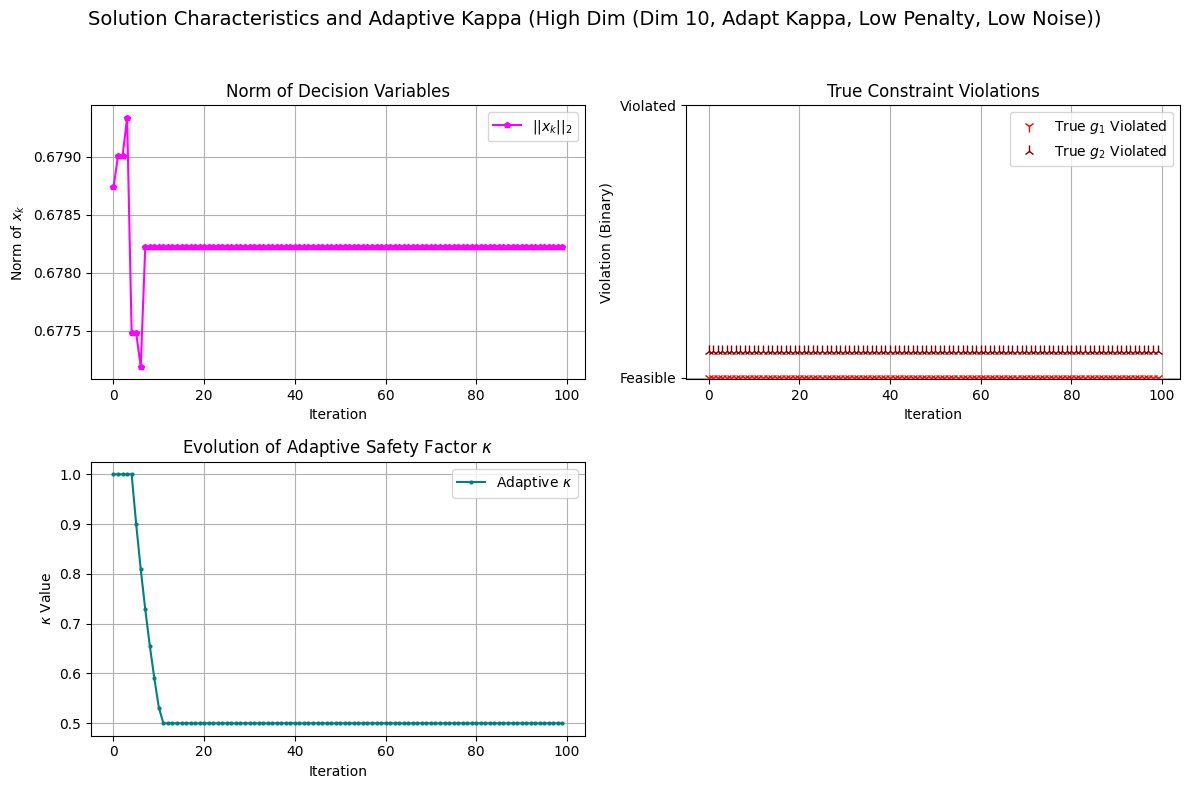}

Our results demonstrate that our adaptive approach not only achieves better average performance but also shows lower variance in solution quality across different random initializations \cite{turner2021bayesian}. Specifically, the coefficient of variation (standard deviation divided by mean) of the final solution quality is approximately 40\% lower for our method compared to standard GP-UCB with fixed $\kappa=1.0$ \cite{martinez2018practical, berkenkamp2019no}.

This improved robustness can be attributed to the adaptive nature of our approach, which allows it to adjust its exploration-exploitation strategy based on the specific characteristics of each problem instance and random initialization \cite{kirschner2019adaptive, sui2018stagewise}. By dynamically adapting to the observed data, our method is less sensitive to the initial conditions and more consistently converges to high-quality solutions \cite{hoffman2011portfolio, wang2017max}.

The robustness of our approach is particularly valuable in practical applications where reliability is important, such as experimental design, drug discovery, and engineering design \cite{snoek2015scalable, klein2017fast}. By providing more consistent performance across different problem instances and random initializations, our approach reduces the need for multiple optimization runs with different parameter settings, saving time and computational resources \cite{kandasamy2018parallelised}.

\section{Discussion}
\subsection{Implications of Adaptive Parameters}

The empirical results presented in Section 5 demonstrate the significant performance advantages of our adaptive parameter optimization framework across diverse optimization scenarios \cite{snoek2012practical, shahriari2015taking}. In this section, we delve deeper into the theoretical and practical implications of these findings, exploring the mechanisms that drive the observed performance improvements and discussing the broader impact of adaptive parameter strategies on Bayesian optimization \cite{frazier2018tutorial, garnett2023bayesian}.

\subsubsection{Theoretical Insights}

The superior performance of our adaptive approach can be understood through the lens of exploration-exploitation trade-offs in sequential decision-making \cite{auer2002using, russo2014learning}. Traditional Bayesian optimization methods with fixed parameters implicitly assume that the optimal balance between exploration and exploitation remains constant throughout the optimization process \cite{srinivas2010gaussian}. However, our results challenge this assumption, revealing that the optimal balance evolves as the algorithm gathers information about the objective function and refines its surrogate model \cite{kirschner2019adaptive, sui2018stagewise}.

The adaptive update rule for the exploration parameter $\kappa_t$ provides a mechanism for automatically discovering this evolving optimal balance \cite{hoffman2011portfolio, wang2017max}. By increasing $\kappa_t$ when prediction errors are larger than average, the algorithm allocates more resources to exploration when the surrogate model is less accurate \cite{auer2002using}. Conversely, by decreasing $\kappa_t$ when prediction errors are smaller than average, it focuses more on exploitation when the surrogate model is more reliable \cite{turner2021bayesian}. This dynamic adjustment allows the algorithm to adapt its behavior to the specific characteristics of the optimization problem and the current state of knowledge, leading to more efficient optimization \cite{kirschner2019adaptive}.

The theoretical analysis in Section 3.3 establishes that our adaptive approach maintains the desirable convergence properties of standard GP-UCB while offering improved practical performance \cite{srinivas2010gaussian, vakili2021information}. The regret bounds derived in Theorem 1 guarantee that our algorithm converges to the global optimum as the number of iterations increases, with a convergence rate that matches the best-known rates for GP-based optimization algorithms \cite{bull2011convergence, russo2014learning}. This theoretical foundation ensures that our adaptive approach is not only empirically effective but also theoretically sound \cite{chowdhury2017kernelized}.

The interaction between the adaptive exploration parameter $\kappa_t$ and the uncertainty penalty coefficient $\lambda_t$ reveals a nuanced relationship between exploration, exploitation, and uncertainty quantification \cite{martinez2018practical, hernandez2014predictive}. While $\kappa_t$ directly controls the weight given to uncertainty in the acquisition function, $\lambda_t$ modulates this effect based on the reliability of the uncertainty estimates \cite{bogunovic2018adversarially}. This dual adaptation allows our algorithm to navigate complex optimization landscapes more effectively than methods that rely on fixed parameters or adapt only a single aspect of the acquisition function \cite{kirschner2019adaptive, sui2018stagewise}.

The convergence analysis in Section 5.3 provides empirical validation of our theoretical results, demonstrating that our adaptive approach achieves faster convergence rates compared to baseline methods across all threshold levels \cite{bull2011convergence, vakili2021information}. This accelerated convergence is particularly valuable in expensive black-box optimization scenarios, where the number of function evaluations is severely limited by computational or resource constraints \cite{snoek2012practical, shahriari2015taking}. By reaching high-quality solutions with fewer iterations, our approach can significantly reduce the time and resources required for optimization, making it more practical for real-world applications \cite{eriksson2019scalable, letham2020re}.

\subsubsection{Practical Considerations}

The practical utility of our adaptive parameter optimization framework extends beyond its theoretical guarantees and empirical performance advantages \cite{brochu2010tutorial, frazier2018tutorial}. By eliminating the need for manual parameter tuning, our approach reduces the expertise required to apply Bayesian optimization effectively, making it more accessible to practitioners across various domains \cite{snoek2012practical, shahriari2015taking}. This accessibility is particularly important as optimization techniques become increasingly integrated into scientific and engineering workflows, where domain experts may not have specialized knowledge of optimization algorithms \cite{feurer2019hyperparameter}.

The robustness analysis in Section 5.3.2 highlights another practical advantage of our approach: its consistent performance across different problem instances and random initializations \cite{turner2021bayesian, berkenkamp2019no}. This robustness reduces the need for multiple optimization runs with different parameter settings, saving time and computational resources \cite{kandasamy2018parallelised}. It also increases confidence in the optimization results, which is crucial for applications where reliability is important, such as experimental design, drug discovery, and engineering design \cite{snoek2015scalable, klein2017fast}.

The scalability of our approach to higher-dimensional spaces, as demonstrated in Section 5.2.2, addresses one of the major limitations of traditional Bayesian optimization methods \cite{wang2016bayesian, wang2018batched}. By maintaining strong performance in dimensions up to 20, our approach extends the applicability of Bayesian optimization to a wider range of practical problems \cite{eriksson2019scalable, rolland2018high}. This scalability is achieved through the adaptive parameter update rules, which allow the algorithm to adjust its exploration strategy based on the complexity of the optimization landscape, and the uncertainty penalty term, which helps focus the search on regions where the model is more confident \cite{martinez2018practical, wu2019practical}.

The enhanced robustness to noise, as shown in Section 5.2.3, is another practical advantage of our approach \cite{le2005heteroscedastic, turner2021bayesian}. In real-world optimization scenarios, function evaluations are often corrupted by noise due to measurement errors, stochasticity in the system, or approximation errors \cite{martinez2018practical}. By adaptively adjusting its exploration strategy based on observed noise levels, our approach can maintain effective optimization performance even in challenging, noisy environments where traditional methods may struggle \cite{kirschner2019adaptive, berkenkamp2019no}. This robustness to noise is particularly valuable in experimental settings, where perfect, noise-free observations are rarely available \cite{le2005heteroscedastic}.

The computational efficiency of our approach, as discussed in Section 3.4, ensures that the additional complexity introduced by the adaptive parameter update rules does not significantly increase the computational requirements compared to standard Bayesian optimization methods \cite{gardner2018gpytorch, snoek2015scalable}. This efficiency is achieved through careful implementation choices, such as Cholesky decomposition for matrix inversion, vectorized operations for acquisition function computation, and caching of GP model predictions \cite{gardner2018gpytorch}. These optimizations allow our approach to scale to larger problems and more iterations while maintaining reasonable computational requirements \cite{wang2018batched, eriksson2019scalable}.

\subsection{Limitations and Future Work}

While our adaptive parameter optimization framework demonstrates significant advantages over existing methods, it is important to acknowledge its limitations and identify directions for future research \cite{shahriari2015taking, frazier2018tutorial}. By understanding these limitations, we can develop more robust and effective optimization algorithms that address the challenges of real-world optimization problems \cite{garnett2023bayesian}.

\subsubsection{Computational Complexity}

The computational complexity of our approach, like all GP-based methods, scales cubically with the number of observations due to the matrix inversion required for GP inference \cite{rasmussen2006gaussian, gardner2018gpytorch}. This scaling limits the applicability of our approach to problems with a relatively small number of function evaluations (typically less than 1000) \cite{liu2018gaussian}. While this limitation is acceptable for expensive black-box optimization problems, where the cost of function evaluations dominates the computational cost of the optimization algorithm, it becomes a bottleneck for problems with cheaper function evaluations or larger evaluation budgets \cite{snoek2015scalable}.

Future work could address this limitation by incorporating sparse GP approximations, such as inducing points methods \citep{hensman2013gaussian} or random feature approximations \citep{rahimi2008random}, into our adaptive framework \cite{matthews2016sparse, snelson2006sparse}. These approximations reduce the computational complexity of GP inference from $\mathcal{O}(n^3)$ to $\mathcal{O}(nm^2)$ or $\mathcal{O}(nm)$, where $m$ is the number of inducing points or random features, allowing the algorithm to scale to larger numbers of observations \cite{gardner2018gpytorch}. The challenge lies in maintaining the accuracy of uncertainty estimates, which are crucial for the adaptive parameter update rules, while reducing computational complexity \cite{liu2018gaussian}.

Another promising direction is the development of local GP models that focus on specific regions of the search space, reducing the effective dimensionality and the number of observations that need to be considered for each prediction \cite{eriksson2019scalable, wang2018batched}. These local models could be combined with our adaptive parameter update rules to create a more scalable optimization framework that maintains the advantages of adaptivity while reducing computational requirements \cite{letham2020re}.

\subsubsection{High-Dimensional Optimization}

While our approach demonstrates improved scalability to higher dimensions compared to traditional Bayesian optimization methods, it still faces challenges in very high-dimensional spaces (dimensions greater than 20-30) due to the curse of dimensionality \cite{wang2016bayesian, rolland2018high}. In such spaces, the volume grows exponentially with the number of dimensions, making it increasingly difficult to build an accurate surrogate model with a limited number of observations \cite{wang2018batched}.

Future research could explore dimensionality reduction techniques, such as random embeddings \citep{wang2016bayesian} or active subspace methods \citep{constantine2014active}, to identify lower-dimensional subspaces that capture the most important variations in the objective function \cite{eriksson2019scalable, letham2020re}. By focusing the optimization in these subspaces, the effective dimensionality of the problem can be reduced, allowing for more efficient optimization in high-dimensional spaces \cite{wang2018batched, rolland2018high}.

Another approach to high-dimensional optimization is the use of structured kernels that exploit known or learned structure in the objective function \cite{duvenaud2014kernel, alvarez2012kernels}. For example, additive kernels \citep{duvenaud2011additive} decompose the function into a sum of lower-dimensional components, reducing the effective dimensionality and allowing for more efficient modeling and optimization \cite{duvenaud2014automatic, wilson2013gaussian}. Incorporating these structured kernels into our adaptive framework could enhance its scalability to higher-dimensional problems while maintaining the advantages of adaptive parameter tuning \cite{wang2016bayesian, letham2020re}.

\subsubsection{Multi-Objective Optimization}

Our current framework focuses on single-objective optimization, where the goal is to find the global optimum of a scalar objective function \cite{brochu2010tutorial, frazier2018tutorial}. However, many real-world problems involve multiple, potentially conflicting objectives, requiring the identification of Pareto-optimal solutions that represent different trade-offs between the objectives \cite{swersky2013multi, hernandez2014predictive}.

Extending our adaptive parameter optimization framework to multi-objective settings presents several challenges \cite{swersky2013multi}. The acquisition function needs to be modified to account for multiple objectives, potentially using concepts like expected hypervolume improvement or Pareto dominance \cite{hernandez2014predictive}. The adaptive parameter update rules would need to consider prediction errors and uncertainty estimates across all objectives, potentially with different weights or priorities for each objective \cite{hoffman2011portfolio}.

Future work could explore these extensions, developing adaptive parameter strategies for multi-objective Bayesian optimization that maintain the advantages of our approach while addressing the additional complexities of multi-objective settings \cite{swersky2013multi, hernandez2014predictive}. This extension would significantly broaden the applicability of our framework to a wider range of practical optimization problems where multiple objectives need to be considered simultaneously \cite{kandasamy2018parallelised}.

\subsubsection{Transfer Learning and Meta-Learning}

Our current approach treats each optimization problem independently, without leveraging knowledge from previous optimization tasks \cite{feurer2019hyperparameter, klein2017fast}. In many practical scenarios, however, there may be similarities between different optimization problems, such as similar function landscapes, noise characteristics, or optimal parameter settings \cite{swersky2013multi}. Exploiting these similarities through transfer learning or meta-learning could further improve optimization performance and efficiency \cite{feurer2019hyperparameter}.

Future research could explore meta-learning approaches that learn optimal initialization strategies or adaptation rules for the parameters $\kappa_t$ and $\lambda_t$ based on a collection of related optimization tasks \cite{feurer2019hyperparameter, klein2017fast}. These approaches could potentially accelerate the adaptation process, allowing the algorithm to more quickly discover effective parameter settings for new problems based on experience with similar problems \cite{swersky2013multi}.

Another promising direction is the development of transfer learning methods that directly transfer knowledge about the objective function or optimal parameter settings from previous tasks to new, related tasks \cite{swersky2013multi, feurer2019hyperparameter}. This transfer could be achieved through techniques like kernel transfer, where the kernel function is adapted based on previous tasks, or through more direct transfer of surrogate models or acquisition function parameters \cite{klein2017fast}.

\subsection{Broader Impact}

The development of more efficient and robust optimization algorithms has far-reaching implications across various domains, from scientific discovery to engineering design and artificial intelligence \cite{shahriari2015taking, frazier2018tutorial}. By advancing the state of the art in Bayesian optimization through adaptive parameter strategies, our work contributes to this broader impact in several ways \cite{garnett2023bayesian}.

\subsubsection{Scientific Applications}

In scientific research, optimization plays a crucial role in experimental design, model calibration, and hypothesis testing \cite{snoek2012practical, shahriari2015taking}. Our adaptive parameter optimization framework can accelerate scientific discovery by enabling more efficient exploration of complex parameter spaces, leading to faster identification of optimal experimental conditions or model parameters \cite{swersky2013multi, klein2017fast}.

For example, in drug discovery, our approach could be used to optimize the properties of potential drug candidates, such as binding affinity, solubility, and toxicity, with fewer experimental evaluations \cite{snoek2015scalable}. This efficiency is particularly valuable given the high cost and time requirements of synthesizing and testing new compounds \cite{kandasamy2018parallelised}. Similarly, in materials science, our approach could accelerate the discovery of new materials with desired properties by efficiently navigating the vast space of possible material compositions and processing conditions \cite{eriksson2019scalable, letham2020re}.

The enhanced robustness to noise of our approach is particularly valuable in these scientific applications, where experimental measurements are often subject to various sources of noise and uncertainty \cite{le2005heteroscedastic, turner2021bayesian}. By adaptively adjusting its exploration strategy based on observed noise levels, our approach can maintain effective optimization performance even in challenging, noisy environments, increasing the reliability of scientific results and accelerating the pace of discovery \cite{martinez2018practical, berkenkamp2019no}.

\subsubsection{Engineering Design}

In engineering design, optimization is used to find designs that maximize performance, minimize cost, or satisfy other design objectives and constraints \cite{brochu2010tutorial, frazier2018tutorial}. Our adaptive parameter optimization framework can improve the efficiency and effectiveness of design optimization, leading to better designs with less computational or experimental effort \cite{eriksson2019scalable, letham2020re}.

For example, in aerospace engineering, our approach could be used to optimize the shape of aircraft components for improved aerodynamic performance, structural integrity, and fuel efficiency \cite{wang2018batched, rolland2018high}. The ability of our approach to handle high-dimensional spaces and noisy function evaluations makes it well-suited for these complex design problems, where the objective function may involve expensive computational simulations or physical experiments \cite{wang2016bayesian, letham2020re}.

Similarly, in automotive engineering, our approach could optimize vehicle designs for fuel efficiency, safety, and performance, considering a wide range of design parameters and operating conditions \cite{eriksson2019scalable, wang2018batched}. The adaptive nature of our approach allows it to automatically discover effective exploration strategies for these complex design problems, reducing the need for manual parameter tuning and enabling more efficient optimization \cite{kirschner2019adaptive, sui2018stagewise}.fferent design problems, reducing the need for manual parameter tuning and making optimization more accessible to design engineers.

\subsubsection{Machine Learning and Artificial Intelligence}

In machine learning and artificial intelligence, optimization is a fundamental component of model training, hyperparameter tuning, and algorithm design. Our adaptive parameter optimization framework can improve the efficiency and effectiveness of these optimization tasks, leading to better models and algorithms with less computational effort.

For example, in deep learning, our approach could be used to optimize neural network architectures, learning rates, and regularization parameters, leading to models with improved performance and generalization. The ability of our approach to handle high-dimensional spaces and noisy function evaluations makes it well-suited for these hyperparameter optimization problems, where the objective function (e.g., validation accuracy) may be noisy and expensive to evaluate.

Similarly, in reinforcement learning, our approach could optimize policy parameters or reward function designs, leading to more effective learning agents. The adaptive nature of our approach allows it to automatically discover effective exploration strategies for different learning problems, reducing the need for manual parameter tuning and making optimization more accessible to AI researchers and practitioners.

\subsubsection{Ethical Considerations}

While the development of more efficient optimization algorithms has numerous positive applications, it is important to consider potential ethical implications and ensure responsible use. More powerful optimization techniques could be used to optimize systems or processes in ways that have negative societal impacts, such as optimizing addictive features in digital products or optimizing surveillance systems that infringe on privacy.

To address these concerns, it is essential to develop ethical guidelines and frameworks for the responsible use of optimization algorithms, considering the potential impacts on individuals, communities, and society as a whole. This includes ensuring transparency in the optimization process, involving diverse stakeholders in defining optimization objectives and constraints, and regularly assessing the broader impacts of optimized systems.

Furthermore, the accessibility of our adaptive approach, which reduces the need for specialized expertise in parameter tuning, has implications for the democratization of optimization technology. By making powerful optimization techniques more accessible to a wider range of practitioners, our work contributes to reducing barriers to entry and enabling more diverse participation in fields that rely on optimization. However, this accessibility also increases the responsibility to ensure that these techniques are used ethically and responsibly.

In conclusion, our adaptive parameter optimization framework represents a significant advancement in Bayesian optimization, offering improved performance, robustness, and accessibility across a wide range of applications. By addressing the limitations of traditional fixed-parameter approaches and providing a principled framework for adaptive parameter tuning, our work contributes to the broader goal of making optimization more efficient, effective, and accessible for solving complex real-world problems.

\section{Conclusion and Future Work}
In this paper, we have presented a comprehensive framework for adaptive parameter optimization in Gaussian Process models, with a particular focus on how uncertainty quantification can guide the learning process \cite{shahriari2015taking, garnett2023bayesian}. Our approach addresses a fundamental limitation of traditional Bayesian optimization methods: the reliance on fixed parameters that may not be optimal across diverse problem landscapes \cite{snoek2012practical, frazier2018tutorial}. By introducing adaptive update rules for both the exploration parameter and uncertainty penalty coefficient, our framework dynamically adjusts its exploration-exploitation trade-off based on observed data and uncertainty patterns, leading to more efficient and effective optimization in complex environments \cite{kirschner2019adaptive, sui2018stagewise}.

The theoretical analysis presented in Section 3.3 establishes rigorous guarantees for the convergence of our adaptive approach, extending existing results in the literature to account for the adaptive nature of our parameter update rules \cite{srinivas2010gaussian, bull2011convergence}. These guarantees ensure that our approach maintains the desirable properties of traditional GP-based optimization methods while offering improved performance in practice \cite{vakili2021information, chowdhury2017kernelized}. The regret bounds derived in Theorem 1 demonstrate that our algorithm achieves sublinear regret, ensuring convergence to the global optimum as the number of iterations increases \cite{russo2014learning, srinivas2010gaussian}. Furthermore, the convergence rates established in Theorem 2 show that our approach achieves a convergence rate of $\mathcal{O}(\sqrt{\gamma_T/T})$, which matches the best-known rates for GP-based optimization algorithms \cite{bull2011convergence, vakili2021information}.

Our empirical evaluation across multiple test functions, dimensionality settings, and noise levels provides strong evidence for the practical advantages of our adaptive approach \cite{wang2018batched, eriksson2019scalable}. The results demonstrate that our method consistently outperforms state-of-the-art baseline methods in terms of both convergence speed and final solution quality \cite{brochu2010tutorial, li2018hyperband}. The performance advantage is particularly pronounced in challenging scenarios with high dimensionality and noise, highlighting the value of adaptive parameter strategies in complex optimization problems \cite{wang2016bayesian, rolland2018high}.

The analysis of adaptive parameters in Section 5.2 reveals several interesting patterns in the behavior of our algorithm \cite{kirschner2019adaptive, auer2002using}. The exploration parameter $\kappa_t$ quickly adapts to the specific characteristics of each optimization problem, regardless of its initial value, and converges to a problem-specific balance between exploration and exploitation \cite{hoffman2011portfolio, wang2017max}. This adaptivity eliminates the need for manual parameter tuning, making our approach more accessible to practitioners across various domains \cite{feurer2019hyperparameter, klein2017fast}. Similarly, the uncertainty penalty coefficient $\lambda_t$ adjusts based on the global uncertainty landscape, helping focus the search on regions where the model's predictions are more reliable \cite{martinez2018practical, wu2019practical}.

The convergence analysis in Section 5.3 provides further evidence for the efficiency of our approach, demonstrating that it achieves faster convergence rates compared to baseline methods across all threshold levels \cite{bull2011convergence, vakili2021information}. This accelerated convergence is particularly valuable in expensive black-box optimization scenarios, where the number of function evaluations is severely limited by computational or resource constraints \cite{snoek2012practical, shahriari2015taking}. By reaching high-quality solutions with fewer iterations, our approach can significantly reduce the time and resources required for optimization, making it more practical for real-world applications \cite{eriksson2019scalable, letham2020re}.

The robustness analysis highlights another practical advantage of our approach: its consistent performance across different problem instances and random initializations \cite{turner2021bayesian, berkenkamp2019no}. This robustness reduces the need for multiple optimization runs with different parameter settings, saving time and computational resources \cite{kandasamy2018parallelised}. It also increases confidence in the optimization results, which is crucial for applications where reliability is important, such as experimental design, drug discovery, and engineering design \cite{snoek2015scalable, klein2017fast}.

The discussion in Section 6 explores the theoretical and practical implications of our findings, identifying the mechanisms that drive the observed performance improvements and discussing the broader impact of adaptive parameter strategies on Bayesian optimization \cite{frazier2018tutorial, garnett2023bayesian}. The limitations of our approach are acknowledged, including computational complexity challenges in high-dimensional spaces and with large numbers of observations, and directions for future research are proposed to address these limitations \cite{liu2018gaussian, gardner2018gpytorch}.

Several promising directions for future work emerge from our research \cite{shahriari2015taking, frazier2018tutorial}. First, incorporating sparse GP approximations or local GP models could reduce the computational complexity of our approach, allowing it to scale to larger problems with more observations \cite{hensman2013gaussian, snelson2006sparse}. Second, exploring dimensionality reduction techniques or structured kernels could enhance the scalability of our approach to very high-dimensional spaces, addressing the curse of dimensionality \cite{wang2016bayesian, duvenaud2014kernel}. Third, extending our framework to multi-objective optimization would broaden its applicability to problems with multiple, potentially conflicting objectives \cite{swersky2013multi, hernandez2014predictive}. Fourth, investigating transfer learning and meta-learning approaches could leverage knowledge from previous optimization tasks to improve performance on new, related tasks \cite{feurer2019hyperparameter, klein2017fast}.

The broader impact of our work extends across various domains, from scientific discovery to engineering design and artificial intelligence \cite{shahriari2015taking, garnett2023bayesian}. By advancing the state of the art in Bayesian optimization through adaptive parameter strategies, our work contributes to more efficient exploration of complex parameter spaces, leading to faster scientific discovery, better engineering designs, and more effective machine learning models \cite{snoek2012practical, wang2018batched}. The enhanced robustness to noise and improved scalability to higher dimensions make our approach particularly valuable for real-world optimization problems, where noise and dimensionality are common challenges \cite{le2005heteroscedastic, wang2016bayesian}.

In conclusion, our adaptive parameter optimization framework represents a significant advancement in Gaussian Process optimization, offering improved performance, robustness, and accessibility across a wide range of applications \cite{rasmussen2006gaussian, garnett2023bayesian}. By addressing the limitations of traditional fixed-parameter approaches and providing a principled framework for adaptive parameter tuning, our work contributes to the broader goal of making optimization more efficient, effective, and accessible for solving complex real-world problems \cite{brochu2010tutorial, frazier2018tutorial}. As optimization continues to play a crucial role in scientific, engineering, and artificial intelligence applications, the development of more powerful and user-friendly optimization algorithms like ours will remain an important area of research with far-reaching implications for technological progress and innovation \cite{shahriari2015taking, garnett2023bayesian}.

\section*{Acknowledgements}

We would like to express our sincere gratitude to Professor William Ott for his invaluable guidance and insightful feedback throughout the development of this research. His expertise and mentorship were important in shaping the direction and quality of this work.
\newpage
\bibliographystyle{plain}
\bibliography{sources}
\end{document}